\documentclass[onefignum,onetabnum]{siamart171218}

\usepackage{lipsum}
\usepackage{amsfonts}
\usepackage{graphicx}
\usepackage{epstopdf}
\usepackage{algorithmic}
\ifpdf
  \DeclareGraphicsExtensions{.eps,.pdf,.png,.jpg}
\else
  \DeclareGraphicsExtensions{.eps}
\fi

\newsiamremark{remark}{Remark}
\newsiamremark{hypothesis}{Hypothesis}
\crefname{hypothesis}{Hypothesis}{Hypotheses}
\newsiamthm{claim}{Claim}

\usepackage{amsopn}

 \usepackage{subfigure}
 
 \usepackage{tikz}
\usepackage{verbatim}
\usetikzlibrary{arrows,backgrounds,snakes,shapes} 

\theoremstyle{plain}
\newtheorem{thm}{Theorem}[section]

\theoremstyle{plainNoItalics}

\newtheorem{example}[thm]{Example}

\graphicspath{{./}{./figure/}}

\ifpdf
\hypersetup{
  pdftitle={An ELDG Article},
  pdfauthor={X. Cai, J.-M. Qiu, and Y. Yang}
}
\fi

\begin{document}


\begin{center}
{\bf

}
\end{center}
 

\title{An Eulerian-Lagrangian discontinuous Galerkin method for transport problems and \\
  its application to nonlinear dynamics%
  \thanks{%
\funding{Research of the first and second author is supported by NSF grant NSF-DMS-1818924, Air Force Office of Scientific Computing FA9550-18-1-0257 and University of Delaware.
Research of the third author is supported by NSF grant DMS-1818467.
 }} }

\author{Xiaofeng Cai
  \thanks{Department of Mathematical Sciences, University of Delaware, Newark, DE, 19716, USA. (\email{xfcai@udel.edu}).}
  \and
  Jing-Mei Qiu%
  \thanks{Department of Mathematical Sciences, University of Delaware, Newark, DE, 19716, USA. (\email{jingqiu@udel.edu}).}
  \and
  Yang Yang%
  \thanks{Department of Mathematical Sciences, Michigan Technological University, Houghton, MI 49931, USA. (\email{yyang7@mtu.edu}).}
}

 \headers{An Eulerian-Lagrangian discontinuous Galerkin method}
{X. Cai, J.-M. Qiu,   and Y. Yang}
\maketitle

 \begin{abstract}
We propose a new Eulerian-Lagrangian (EL) discontinuous Galerkin (DG) method. The method is designed as a generalization of the semi-Lagrangian (SL) DG method for linear advection problems proposed in [J. Sci. Comput. 73: 514-542, 2017], which is formulated based on an adjoint problem and tracing upstream cells by tracking characteristics curves highly accurately. In the SLDG method, depending on the velocity field, upstream cells could be of arbitrary shape. Thus, a more sophisticated approximation to sides of the upstream cells is required to get high order approximation. For example, quadratic-curved (QC) quadrilaterals were proposed to approximate upstream cells for a third-order spatial accuracy in a swirling deformation example. In this paper, for linear advection problems, we propose a more general formulation, named the ELDG method. The scheme is formulated based on a {\em modified} adjoint problem for which the upstream cells are always quadrilaterals, which avoids the need to use QC quadrilaterals in the SLDG algorithm. The newly proposed ELDG method can be viewed as a new general framework, in which both the classical Eulerian Runge-Kutta DG formulation and the SL DG formulation can fit in. Numerical results on linear transport problems, as well as the nonlinear Vlasov and incompressible Euler dynamics using the exponential RK time integrators, are presented to demonstrate the effectiveness of the ELDG method.
 \end{abstract}

\begin{keywords}
 Eulerian-Lagrangian; discontinuous Galerkin; mass conservative; semi-Lagrangian; Vlasov simulations;  characteristics.
\end{keywords}

\begin{AMS}
  65M25, 65M60, 76M10
  
\end{AMS}

\section{Introduction}

We propose a new Eulerian-Lagrangian (EL) discontinuous Galerkin (DG) method for a model transport equation in the form of
 \begin{equation}
u_t + \nabla \cdot( \mathbf{P}( u;\mathbf{x},t ) u ) = 0, \ (\mathbf{x},t)\in\mathbb{R}^d\times[0,T],
\label{general_nonlinear}
\end{equation}
which could come from a wide range of application fields including fluid dynamics, climate modeling, and kinetic description of plasma.
There are three main classes of computational methods for solving \eqref{general_nonlinear}: Lagrangian, Eulerian and semi-Lagrangian (SL). Each class of methods has their own advantages and limitations. The Lagrangian method is particle based, works efficiently for high dimensional problems, but suffers from statistical noises; while the latter two methods are mesh-based method, can be designed to be of high order accurate, but suffers from the curse of dimensionality. The main difference between Eulerian and SL methods is the space-time region in consideration: the Eulerian method performs numerical discretizations with fixed spatial locations in time; while the semi-Lagrangian method usually do that along convection characteristics. When characteristics are tracked accurately, semi-Lagrangian methods often allow much larger time stepping sizes than their Eulerian counterparts.

Among different classes of SL methods in the literature, we would like to mention a few closely related ones that are developed in the finite element framework. There is a line of research work along Eulerian Lagrangian Localized Adjoint Methods (ELLAM) \cite{celia1990eulerian}. ELLAM introduces an adjoint problem for the test function in the continuous finite element framework and has a broad range of influence in different application fronts \cite{wang1999family, russell2002overview}. Compared with ELLAM, the SLDG \cite{cai2016high} is being developed in the discontinuous Galerkin finite element framework.  SL schemes could be developed base on forward \cite{bosler2019conservative} or backward characteristics tracing.
  Here we choose to develop our schemes base on backward characteristics tracing.

In this paper, we propose a new ELDG method that is mesh-based, and is a generalized framework of the SL DG method developed earlier \cite{cai2016high}. It is designed to take advantage of information propagation along characteristics as in a SL method, and maintain essential properties of the SLDG method on mass conservation, high order spatial and temporal accuracy, and allowing for extra large time steps with stability.
We first focus on developing the ELDG algorithm for linear transport problems. A new ingredient of the method is the introduction of a {\em modified} adjoint problem for the test function. {\em The velocity field of the modified adjoint problem is a linear function that approximates that of the original transport problem.} There are two positive consequences of such modification. One is that the test function remains in the same $P^k$ polynomial spaces, whereas in the SLDG setting the test function does not necessarily remain in $P^k$ and needs to be approximated. In fact, a close connection can be drawn between the ELDG method and the Arbitrary Lagrangian Eulerian (ALE) DG method \cite{klingenberg2017arbitrary}, when we view the space-time region in the ELDG method as a dynamic moving mesh.
The second advantage brought by the modified adjoint problem is that the shape of upstream cells is always quadrilaterals in a 2D setting. For a general variable coefficient problem, upstream cells of the SLDG method could be of arbitrary shape and needs to be better approximated. In \cite{cai2016high}, we propose to use quadratic curves in approximating sides of upstream cells, so that we have third order spatial accuracy. Such a practice is difficult be further generalized to schemes with even higher order accuracy, and for problems in higher-dimensions. With the newly ELDG method, no curves are needed to better approximate upstream cells. A direct generalization of the algorithm to higher dimensional problems can be similarly done in principle.

Due to the approximate nature of the velocity field in the modified adjoint problem, there is an extra flux term taking account of the difference between velocity fields from the modified adjoint problem and the original problem. The newly proposed ELDG scheme evolves this extra flux term in a similar spirit to the classical Eulerian RKDG method \cite{cockburn1989tvb}. The ELDG scheme is designed base on the integral form of the equation over characteristics-related space-time regions; yet we transform such integral formulation into a time-differential form, for which the method-of-lines strong-stability preserving (SSP) Runge-Kutta (RK) can be directly applied. Here, we would like to mention the Eulerian Lagrangian weighted essentially non-oscillatory schemes developed in \cite{huang2017eulerian,huang2018implicit}, for which a different way of treating time integration is proposed.

As nonlinear applications of the ELDG algorithm, we consider the nonlinear Vlasov-Poisson system, the guiding center Vlasov model as well as the incompressible Euler equations. Here, we couple the ELDG algorithm with the RK exponential integrator \cite{celledoni2009semi, cai2019exp} to realize a uniformly high order spatial-temporal discretization of nonlinear transport. In particular, the RK exponential integrator decomposes a time step evolution of the nonlinear problem into the composition of a sequence of linear problems. Extensive numerical experiments are performed and effectiveness of the ELDG method is showcased in various settings with allowance of extra large time stepping sizes.

This paper is organized as follows. In Section \ref{section:1d},   we derive the formulation of ELDG for one-dimensional (1D) linear transport problems,  where the main spirit of the method is introduced.
In Section \ref{section:2d}, we perform a nontrivial generalization of  the scheme for 2D linear transport problems.
In Section \ref{section:nonlinear}, we combine the ELDG scheme with the Runge-Kutta exponential integrators for   nonlinear transport problems.
In Section \ref{section:numerical}, the performance of the proposed method is shown through extensive
numerical tests. Finally, concluding remarks are made in Section \ref{section:conclusion}.

\section{ELDG formulation for 1D linear transport problems} \label{section:1d}

To illustrate the key idea of the ELDG scheme, we start from a 1D linear transport equation in the following form
\begin{equation}
u_t+(a(x, t)u)_x = 0, \quad x\in[x_a, x_b].
\label{scalar1d}
\end{equation}
For simplicity, we assume periodic boundary conditions, and the velocity field $a(x, t)$ is a continuous function of space and time. We perform a partition of the computational domain $x_a=x_{\frac12}< x_{\frac32}<\cdots< x_{N+\frac12} =x_b$.
Let $I_j=[x_{j-\frac12}, x_{j+\frac12} ]$ denote an element of length $\Delta x_j=x_{j+\frac12}-x_{j-\frac12}$ and define $\Delta x=\max_{j}\Delta x_j.$
We define the finite dimensional approximation space, $V_h^k = \{ v_h:  v_h|_{I_j} \in P^k(I_j) \}$, where $P^k(I_j)$ denotes the set of polynomials of degree at most $k$.
For this finite-dimensional space, we introduce a set of basis functions $\{ \Psi_{j,m}(x) \}_{ 1\leq j\leq N, 0\leq m\leq k }$.
 We also introduce a set of basis functions $\{ \psi_{j,m}(x,t) \}_{ 1\leq j\leq N, 0\leq m\leq k }$, which will be used in an adjoint problem.
The subscripts of $\Psi_{j,m}(x) $ and $\psi_{j,m}(x,t) $  are often omitted, when there is no risk of ambiguity.
Moreover, we define $t^n$ to be the $n-$th time level, and $\Delta t =t^{n+1}-t^n$ to be the time-stepping size.

\subsection{Review of SLDG scheme \cite{cai2016high}.}

\begin{figure}[h!]
\centering
\begin{tikzpicture}[x=0.9cm,y=1cm]
  \begin{scope}[thick]

  \draw[white,fill=blue!4] (0.,2) to[out=240,in=85] (-1.5+0.7,0) -- (1.6+0.4,0) to[out=85,in=245] (2.4,2)
      -- cycle;
   \draw (-3,3) node[fill=white] {};
    \draw (-3,-1) node[fill=white] {};
    \draw[black]                   (-3,0) node[left] {} -- (3,0)
                                        node[right]{$t^{n}$};
    \draw[black] (-3,2) node[left] {$$} -- (3,2)
                                        node[right]{$t^{n+1}$};
     \draw[snake=ticks,segment length=2.4cm] (-2.4,2) -- (0,2) node[above left] {$x_{j-\frac12}$};
     \draw[snake=ticks,segment length=2.4cm] (0,2) -- (2.4,2) node[above right] {$x_{j+\frac12}$};

          \draw[snake=ticks,segment length=2.4cm] (-2.4,0) -- (0,0);
     \draw[snake=ticks,segment length=2.4cm] (0,0) -- (2.4,0);
            \draw[-latex,dashed,blue] (0.,2) node[left] {$$} to[out=240,in=85] (-1.5+0.7,0)
                                        node[black,below]{$x_{j-\frac12}^\star$ };

                \draw[-latex,dashed,blue] (2.4,2) node[left] {$$} to[out=245,in=85] (1.6+0.4,0)
                                        node[black,below right]{ $x_{j+\frac12}^\star$};
\draw [decorate,color=red,decoration={brace,mirror,amplitude=9pt},xshift=0pt,yshift=0pt]
(-1.5+0.7,0) -- (1.6+0.4,0) node [red,midway,xshift=0cm,yshift=-20pt]
{\footnotesize $\widetilde{I}_j(t^{n}) = I_{j}^\star$};

\draw [decorate,color=red,decoration={brace,amplitude=10pt},xshift=0pt,yshift=0pt]
(0,2) -- (2.4,2) node [red,midway,xshift=0cm,yshift=20pt]
{\footnotesize $\widetilde{I}_j(t^{n+1}) = I_{j}$};
\draw[-latex,thick,blue](0.2,1)node[right,scale=1.]{$\widetilde{I}_j(t)$}
        to[out=180,in=0] (-0.55 ,1) node[above left=2pt] {$$};  

\draw[-latex,thick,blue](1.1,1)node[right,scale=1.]{ }
        to[out=0,in=180] (2.1,1) node[above left=2pt] {$$};  

  \node[blue, rotate=0] (a) at (-2.  ,1.3)  { \tiny $ \frac{d\tilde{x}(t ) }{dt}= a(\tilde{x}(t ),t)$ };
  \node[blue, rotate=0] (a) at (-2  ,1.3-0.5)  {\tiny where $\tilde{x}(t^{n+1} )={x_{j-\frac12}}$ };

  \draw[-latex,blue ]( -1, 1.3 )node[right=-3pt,scale=1.0]{ }
        to[out=340,in=160] ( -0.5, 1.1) node[below=2pt] { };

  \end{scope}
\end{tikzpicture}
\label{schematic_1d}
\caption{Illustration for the space-time region for the SLDG formulation.  }
\end{figure}
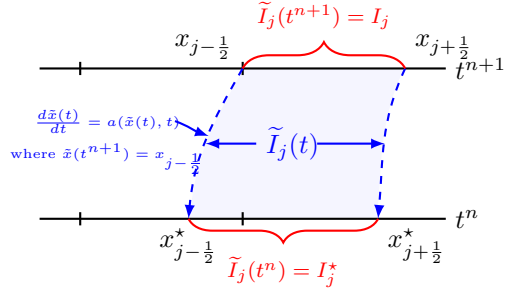

The SLDG method proposed in \cite{cai2016high} is formulated based on an adjoint problem of \eqref{scalar1d} with $\forall \Psi \in P^k(I_j)$,
\begin{equation}
\begin{cases}
\psi_t + a(x,t) \psi_x =0 ,\ t\in[t^n,t^{n+1} ],\\
\psi(t=t^{n+1}) = \Psi(x),
\end{cases}
\label{final-value}
\end{equation}
for which the solution $\psi$ stays constant along characteristic trajectories. It was shown in \cite{Guo2013discontinuous} that
\begin{equation}
\frac{d}{dt}\int_{\widetilde{I}_j(t)} u(x,t)\psi(x,t) dx =0,
\label{property}
\end{equation}
where $\widetilde{I}_j(t)$ is a dynamic interval bounded by characteristics emanating from cell boundaries of $I_j$ at $t=t^{n+1}$, see Figure~\ref{schematic_1d} for illustration. An SL time discretization of \eqref{property} leads to
\begin{equation}
\int_{I_j} u^{n+1} \Psi dx = \int_{I_{j}^\star } u(x,t^n) \psi(x,t^n)dx,
\label{eq: integral1}
\end{equation}
where $I_{j}^\star  = [x_{j-\frac12}^\star , x_{j+\frac12}^\star]$ with $x_{j\pm\frac12}^\star=\widetilde{x}_{j\pm\frac12}(t^n)$ being the foots of trajectory at $t^n$ emanating from
$(x_{j\pm\frac12} ,t^{n+1})$. In order to update the numerical solution $u^{n+1}$, we vary the test function $\Psi$ as basis of $V_h^k$ and evaluate the right-hand side (RHS) integral of \eqref{eq: integral1} properly. The detailed procedures can be found in \cite{cai2016high}.

\subsection{The new ELDG scheme.}

 The newly proposed ELDG method differs from the SLDG method \cite{cai2016high} in the formulation of a {\em modified} adjoint problem for the test function $\psi$. To introduce the scheme, we first introduce the modified adjoint problem and the associated space-time region $\Omega_j$; then we derive a semi-discrete version of the ELDG scheme based on the space-time region of $\Omega_j$; finally a method-of-times Runge-Kutta method is applied for time marching.

\noindent
 {\bf (1) A modified adjoint problem.} We consider the adjoint problem with $\forall \Psi \in P^k(I_j)$
    on the time interval $[t^n, t^{n+1}]$:
    \begin{equation}
\begin{cases}
\psi_t + \alpha(x,t) \psi_x =0 ,\ t\in[t^n,t^{n+1} ],\\
\psi(t=t^{n+1}) = \Psi(x),
\end{cases}
\label{adjoint1d}
\end{equation}
with $\alpha(x,t)$ being a bilinear function of $(x, t)$ designed by three steps below:
\begin{enumerate}
\item {\bf On $I_j$ at $t^{n+1}$:} we let $\alpha(x, t^{n+1})$ be a linear polynomial on $I_j$ interpolating $a(x, t^{n+1})$ at cell boundaries,
\begin{equation}
\label{eq: a}
\alpha(x_{j\pm\frac12},t^{n+1}) =a(x_{j\pm\frac12}, t^{n+1})\doteq \nu_{j\pm\frac12}.
\end{equation}
That is,
 \begin{equation}
 \alpha(x, t^{n+1})=-\nu_{j-\frac12} \frac{x-x_{j+\frac12}}{\Delta x_{j}} + \nu_{j+\frac12} \frac{x-x_{j-\frac12}}{\Delta x_{j}}\in P^1(I_j).
 \label{eq: alpha_def1}
 \end{equation}
 \item We define a space-time region $\Omega_j = \tilde{I}_j(t) \times [t^n, t^{n+1}]$ with  the dynamic interval, $\tilde{I}_j(t)   = [\tilde{x}_{j-\frac12}(t), \tilde{x}_{j+\frac12}(t)]$,$t \in [t^n, t^{n+1}]$, where $\tilde{x}_{j\pm\frac12}(t) = x_{j\pm\frac12} + (t-t^{n+1}) a(x_{j\pm\frac12}, t^{n+1})$ emanating from cell boundaries $x_{j\pm\frac12}$ with slopes $a(x_{j\pm\frac12}, t^{n+1})$.
It will become clear after the third step that the space-time region $\Omega_j$ is the dynamic characteristic region of the modified adjoint problem \eqref{adjoint1d}. We let ${I}^\star_j \doteq \tilde{I}_j(t^n)$ be the upstream cell of $I_j$ at $t^n$.
See the left panel in Figure~\ref{iso} for illustration.
\item {\bf On $\tilde{I}_j(t)$ for $[t^{n}, t^{n+1})$:} let $\tilde{x}(t;(\xi,t^{n+1}) )$ be a straight line emanating from any point $\xi\in I_j$ at $t^{n+1}$ and with the slope $\alpha(\xi, t^{n+1})$.
That is,
\begin{equation}
\label{eq: char_adjoint}
\frac{d}{d t}\tilde{x}(t;(\xi,t^{n+1}) ) =\alpha(\xi, t^{n+1}), \quad \tilde{x}(t^{n+1};(\xi,t^{n+1}) ) = \xi.
\end{equation}
Then
\begin{equation}
\label{eq: tilde_x_1d}
\tilde{x}(\tau;(\xi,t^{n+1})) = \xi - \alpha(\xi, t^{n+1}) (t^{n+1}-\tau), \quad \forall \tau\in[t^n, t^{n+1}).
\end{equation}
We let
\begin{equation}
\alpha(\tilde{x}(\tau;(\xi,t^{n+1})), \tau) = \alpha(\xi, t^{n+1}), \quad \tau\in[t^n, t^{n+1}).
\label{eq: alpha_def2}
\end{equation}
\end{enumerate}
We would like to point out a few facts about $\Omega_j$ and the modified adjoint problem \eqref{adjoint1d}:
\begin{itemize}
\item From the construction of $\Omega_j$ and $\alpha(x, t)$ of the modified adjoint problem \eqref{adjoint1d}, it can be easily checked that, \eqref{eq: char_adjoint} is the characteristics equation for the modified adjoint problem \eqref{adjoint1d}.
\item
$\tilde{x}(\tau;(\xi,t^{n+1}))$ satisfying eq.~\eqref{eq: tilde_x_1d} is a linear function of $\xi$ and $\tau$; the Jacobian is
\begin{equation}
\frac{\partial \tilde{x}(\tau;(\xi,t^{n+1}))}{\partial\xi} = 1-\frac{\nu_{j+\frac12}-\nu_{j-\frac12}}{\Delta x_{j}}(t^{n+1}-\tau),
\label{eq: 1d_jacobian}
\end{equation}
which will become useful later in implementation. In particular,
\begin{equation*}
\frac{\partial \tilde{x}(t^n;(\xi,t^{n+1}))  }{\partial\xi} = 1-\Delta t\frac{\nu_{j+\frac12}-\nu_{j-\frac12}}{\Delta x_j}.
\end{equation*}
\item In order for the characteristics not crossing each other, one has to enforce the condition of $\frac{\partial \tilde{x}(t^n;(\xi,t^{n+1}))  }{\partial \xi} \ge0$, which implies the time step constraint
\begin{equation}
\label{eq: time_constraint}
\Delta t\le \frac{\min_j \Delta x_j}{  \max(\nu_{j+\frac12}-\nu_{j-\frac12},0 ) }.
\end{equation}

\item
For the modified adjoint problem, the solution $\psi$ stays constant along characteristics \eqref{adjoint1d}, therefore we have
\begin{equation}
\psi(\tilde{x}(\tau;  (\xi, t^{n+1})),\tau) = \Psi(\xi) \in P^k(I_j), \quad \forall \tau\in[t^n, t^{n+1}].
\label{eq: test_time}
\end{equation}
If we consider a transformation between $x \in \tilde{I}_j$ to a reference interval $\xi \in I_j$, see Figure~\ref{iso}, eq.~\eqref{eq: test_time} indicates that the test function $\psi(\tilde{x}(\tau;  (\xi, t^{n+1})),\tau)$ in the $\xi$ coordinate remains the same as the classical test function $\Psi(\xi)$, i.e. standard basis functions in $P^k(I_j)$.
\end{itemize}

  \begin{figure}[h!]
\centering
\begin{tikzpicture}[x=0.9cm,y=1cm]
  \begin{scope}[thick]

  \draw[fill=blue!5] (0.,2) -- (-1.5,0) -- (1.6,0) -- (2.4,2)
      -- cycle;
   \node[blue!70, rotate=0] (a) at ( 1. ,1.5) {\LARGE $\Omega_j$ };

   \draw (-3,3) node[fill=white] {};
    \draw (-3,-1) node[fill=white] {};
    \draw[black]                   (-3,0) node[left] {} -- (3,0)
                                        node[right]{$t^{n}$};
    \draw[black] (-3,2) node[left] {$$} -- (3,2)
                                        node[right]{$t^{n+1}$};

     \draw[snake=ticks,segment length=2.4cm] (-2.4,2) -- (0,2) node[above left] {$x_{j-\frac12}$};
     \draw[snake=ticks,segment length=2.4cm] (0,2) -- (2.4,2) node[above right] {$x_{j+\frac12}$};

          \draw[snake=ticks,segment length=2.4cm] (-2.4,0) -- (0,0);
     \draw[snake=ticks,segment length=2.4cm] (0,0) -- (2.4,0);

            \draw[blue,thick] (0.,2) node[left] {$$} -- (-1.5,0)
                                        node[black,below]{$x_{j-\frac12}^\star$ };

                \draw[blue,thick] (2.4,2) node[left] {$$} -- (1.6,0)
                                        node[black,below right]{ $x_{j+\frac12}^\star$};
\draw [decorate,color=red,decoration={brace,mirror,amplitude=9pt},xshift=0pt,yshift=0pt]
(-1.5,0) -- (1.6,0) node [red,midway,xshift=0cm,yshift=-20pt]
{\footnotesize $\widetilde{I}_j(t^{n}) = I_{j}^\star$};

\draw [decorate,color=red,decoration={brace,amplitude=10pt},xshift=0pt,yshift=0pt]
(0,2) -- (2.4,2) node [red,midway,xshift=0cm,yshift=20pt]
{\footnotesize $\widetilde{I}_j(t^{n+1}) = I_{j}$};

\draw[-latex,thick](0,1)node[right,scale=1.]{$\widetilde{I}_j(t)$}
        to[out=180,in=0] (-0.75,1) node[above left=2pt] {$$};  

\draw[-latex,thick](1,1)node[right,scale=1.]{ }
        to[out=0,in=180] (2,1) node[above left=2pt] {$$};  

  \node[blue, rotate=54] (a) at (-1.25,1.) { $\nu_{j-\frac12}$ };
  \node[blue, rotate=70] (a) at (2.35,1.) { $\nu_{j+\frac12}$ };

  \end{scope}

   \draw[-latex,dashed](1.6,1)node[right,scale=1.0]{ }
        to[out=30,in=150] (6,1) node[] {};

     \draw[-latex,red,dashed](1.6,2)node[right,scale=1.0]{ }
        to[out=40,in=140] (6,1.2) node[] {};

     \draw[-latex,red,dashed](1.3,0)node[right,scale=1.0]{ }
        to[out=-60,in=220] (6,1-0.2) node[] {};

\draw[|-|,black,thick] (6.5,0.8) node[below]{\footnotesize $\xi=x_{j-1/2}$} to (9.5,0.8)node[below]{\footnotesize $\xi=x_{j+1/2}$};
\end{tikzpicture}
\label{iso}
\caption{Illustration for the mapping between dynamic element $\widetilde{I}_j(t)$ (left) and the iso-parametric element (right).}
\end{figure}

\bigskip
\noindent{\bf (2) Formulation of the semi-discrete ELDG scheme.} In order to formulate the scheme, we integrate $\text{\eqref{scalar1d}}\cdot\psi + \text{\eqref{adjoint1d}}\cdot u$ over $\Omega_j$, which gives the following identity,
\begin{equation}
\int_{\Omega_j} \left[\text{\eqref{scalar1d}}\cdot\psi + \text{\eqref{adjoint1d}}\cdot u \right] dxdt=0.
\end{equation}
That is,
\begin{align}
0 =& \int_{t^n}^{t^{n+1} } \int_{\tilde{I}_j(t)} \left( u_t\psi+u\psi_t \right)dxdt
+ \int_{t^n}^{t^{n+1} } \int_{\tilde{I}_j(t)}  \left( (a(x, t)u)_x\psi+\alpha(x,t)\psi_x u \right) dxdt \nonumber \\
=&\int_{t^n}^{t^{n+1} } \int_{\tilde{I}_j(t)} ( u\psi )_t dxdt
 + \int_{t^n}^{ t^{n+1} } \int_{ \tilde{I}_j(t) } \left(  (au\psi)_x-au\psi_x + \alpha\psi_x u  \right)dxdt
 \nonumber \\
 =& \int_{t^n}^{t^{n+1} } \left[
\frac{d}{dt} \int_{\tilde{I}_j(t) }u\psi dx - \alpha u\psi\left|_{\tilde{x}_{j-\frac12}(t)}^{\tilde{x}_{j+\frac12}(t)} \right.
+   a u\psi\left|^{\tilde{x}_{j+\frac12}(t) }_{\tilde{x}_{j-\frac12}(t)} \right.
+\int_{ \tilde{I}_j(t) }  (\alpha-a)u \psi_x dx \right]dt
\nonumber\\
=&\int_{t^n}^{t^{n+1}} \left[
\frac{d}{dt} \int_{\tilde{I}_j(t) }u\psi dx + (a-\alpha) u\psi\left|_{\tilde{x}_{j-\frac12}(t)}^{\tilde{x}_{j+\frac12}(t)} \right.
-\int_{\tilde{I}_j(t) }  (a-\alpha)u \psi_x dx
\right]dt.
\label{eq: int}
\end{align}
Letting $F(u) \doteq (a-\alpha)  u$, the time differential form of \eqref{eq: int} gives
\begin{equation}
\frac{d}{dt} \int_{\tilde{I}_j(t)}(u\psi)dx =-  \left(F\psi\right) \left|_{\tilde{x}_{j+\frac12}(t) } \right. +   \left(F\psi\right) \left|_{\tilde{x}_{j-\frac12}(t) } \right.   + \int_{\tilde{I}_j(t)}F\psi_xdx.
\label{mol_2}
\end{equation}
Notice that the dynamic interval of $\tilde{I}_j(t)$ can always be linearly mapped to a reference cell $\xi\in I_j$, see the right plot in Figure~\ref{iso}, then eq.~\eqref{mol_2} in the $\xi$-coordinate becomes
\begin{equation}
\frac{d}{dt} \int_{I_j}(u  \Psi(\xi))\frac{\partial \tilde{x}(t;(\xi,t^{n+1})) }{\partial \xi}d\xi =-  \left(F\Psi\right) \left|_{ \xi=x_{j+\frac12}  } \right. +   \left(F\Psi\right) \left|_{ \xi=x_{j-\frac12}  } \right.      + \int_{I_j}F \Psi_{\xi}d\xi.
\label{mol_4_1}
\end{equation}
The DG discretization \cite{cockburn1991runge,cockburn1989tvb} of \eqref{mol_4_1} is to find $u_h(\xi,t ) \in P^{k}( I_j )$ as the approximate solution of $u(\tilde{x}(t;(\xi,t^{n+1})),t)$ on $\tilde{I}_j(t)$, so that for $\forall \Psi \in P^k(I_j)$,
\begin{equation}
\frac{d}{dt} \int_{I_j}u_h  \Psi \frac{\partial \tilde{x}(t;(\xi,t^{n+1})) }{\partial \xi}d\xi =-  \hat{F}_{j+\frac12}\Psi(x_{j+\frac12}^-) +   \hat{F}_{j-\frac12}\Psi(x_{j-\frac12}^+)     + \int_{I_j}F \Psi_{\xi}d\xi.
\label{mol_4}
\end{equation}
Notice here ${u}_h$ could be discontinuous across $x_{j-\frac12}^\star$.
In this paper, we choose $\hat{F}$ as a monotone flux, e.g. the Lax-Friedrichs flux
\begin{equation}\label{flux}
\hat{F}(u^-, u^+) = \frac12 (F(u^-) + F(u^+))-\frac{\alpha_0}{2}(u^+ -u^-), \quad \alpha_0 = \max_{u} |F'(u)|;
\end{equation}
and we use Gauss quadrature rules with $k+1$ quadrature points to approximate the integral term $\int_{ I_j  }F(u_h)\Psi_\xi d\xi$ on the RHS of the equation \eqref{mol_4}.
%



\bigskip
\noindent{\bf (3) RK time discretization and fully discrete scheme.}
We can write the semi-discrete scheme \eqref{mol_4} into a form of  ordinary differential equations (ODEs) with an initial condition.
We let $\tilde{ \mathbf{U} } (t)$ be a vector in $R^{N(k+1)}$ which consists of degrees of freedom $\{ \int_{ \tilde{I}_j(t) } u_h(x,t) \psi_{j,m}(x,t) dx \doteq \tilde{ U }_{j,m}(t) \}_{1 \leq j\leq N,0 \leq m\leq k}$,
and denote the spatial discretization operator of the RHS of \eqref{mol_4} as $\mathcal{L}\left(\tilde{ \mathbf{U} } (t) ,t \right)$.
Then the semi-discrete scheme \eqref{mol_4} can be written  as
\begin{align}
&\frac{\partial}{\partial t}\tilde{ \mathbf{U}  } (t) = \mathcal{L}\left(\tilde{  \mathbf{U}   } (t) ,t \right), \qquad
\tilde{ \mathbf{U} } (t^n) = \tilde{  \mathbf{U}  }^n.
\label{semi}
\end{align}
There are two main steps involved here.
\begin{enumerate}
\item {\bf Obtain the initial condition} of \eqref{mol_4} by an $L^2$ projection of $u_h$ on upstream cells $\tilde{I}_j$ by SLDG method. In particular, $\tilde{  \mathbf{U}  }^n$ consists of the numerical solutions $ \tilde{U}_{j,m}^n$ of the SLDG scheme \cite{cai2016high} for approximating $$ \int_{ \tilde{I}_j(t^n) } u_h(x,t^n)\psi_{j,m}(x,t^n) dx.$$
\item {\bf Update \eqref{semi} from $\tilde{  \mathbf{U}  }^n$ to $\tilde{  \mathbf{U} }^{n+1}$.} we apply the SSP explicit RK methods \cite{shu1988efficient} as in a method-of-lines approach.
In particular, the time-marching algorithm using an $s$-stage RK method follows the procedure below:
\begin{enumerate}
  \item Get the mesh information of the dynamic element $\tilde{I}_j^{(l)},l=0,\cdots, s$ on RK stages by eq.~\eqref{eq: tilde_x_1d}.
  \item  For RK stages $i=1,\cdots, s$, compute
  \begin{align}
    \tilde{ \mathbf{U} }^{(i)}= \sum_{l=0}^{i-1} \left[ \alpha_{il}   \tilde{  \mathbf{U} }^{(l)}
  + \beta_{il} \Delta t^n \mathcal{L}\left(  \tilde{  \mathbf{U}  }^{(l)},t^n+d_l\Delta t^n \right) \right],
  \label{eq: RK}
  \end{align}
  where $\alpha_{il}$ and $\beta_{il}$ are related to RK methods. They are provided in Table \ref{table:ssprk} for the second order and third order SSP RK methods.
\end{enumerate}
 Note that $\tilde{ \mathbf{U} }^n$ is evaluated by the SLDG scheme in $x$-coordinate, while $\tilde{ \mathbf{U} }^{ (i) }$ in the each time stage is updated with respect to the reference $\xi$ coordinate.
\end{enumerate}
\begin{table}[!ht]
\caption{Parameters of some practical Runge-Kutta time discretizations.  }
\vspace{0.1in}
\centering
\begin{tabular}{ c  cc c }
\hline
Order  &   $\alpha_{il}$  & $\beta_{il}$  & $d_l$    \\
 \hline

 2  &    1  &  1 &  0\\
    & $\frac12$ \ $\frac12$ & 0 \ $\frac12$  &1 \\
\hline
 3  &    1  &  1 &  0\\
    & $\frac34$ \ $\frac14$ & 0 \ $\frac14$  &1 \\
    & $\frac13$ \ 0 \ $\frac23$  & 0 \ 0 \ $\frac23$ & $\frac12$ \\
\hline
\end{tabular}
\label{table:ssprk}
\end{table}

\begin{theorem} (Mass conservation)
Given a DG solution $u_h(x,t^n)\in V_h^k$ and assuming the boundary condition is periodic, the proposed fully discrete ELDG scheme with SSP RK time discretization of \eqref{semi} is locally mass conservative.
In particular,
\begin{equation*}
\sum_{ i=1 }^N \int_{ I_j } u_h( x,t^{n+1} ) dx
=
\sum_{ i=1 }^N \int_{ I_j } u_h( x,t^n ) dx.
\end{equation*}
\end{theorem}
\noindent
{\em Proof.} It can be proved by letting $\psi=1$, the conservative form of integrating $F$ function with unique flux at cell boundaries, as the mass conservation property of SLDG scheme \cite{cai2016high}. We skip details for brevity.

A few remarks are in order for the proposed ELDG scheme, in comparison with existing SLDG  \cite{cai2016high}, RKDG \cite{cockburn1989tvb} and ALE DG \cite{klingenberg2017arbitrary} methods in the literature. These remarks also apply to the 2D ELDG scheme in the next section.

\begin{remark} (Comparison with the SLDG method \cite{cai2016high})
The modified adjoint problem \eqref{adjoint1d} is different from the adjoint problem \eqref{final-value} in the velocity field. In some sense,  $\alpha(x, t)$ is an approximation of $a(x, t)$. While the characteristics induced by $a(x, t)$ could be curves and the test function $\phi$ satisfying eq.~\eqref{final-value} may no longer be polynomials, the characteristics induced by $\alpha(x, t)$ are straight lines and the test function $\phi$ remains a $P^{k}$ polynomial on $\tilde{I}_j(t)$. The difference, between $\alpha(x, t)$ and exact slopes $a(x, t)$ for characteristic curves, is taken into account by the $F$ function in \eqref{mol_4}.
\end{remark}

\begin{remark} (A framework encompassing RKDG and SLDG)
The new scheme formulation \eqref{mol_4} offers a general framework that encompasses the traditional Eulerian RKDG scheme \cite{cockburn2001runge} and the SLDG method proposed in \cite{cai2016high}. For the linear equation with the special case of $\alpha=a$, the ELDG method becomes the SLDG method \cite{cai2016high} and the scheme is unconditionally stable.
 In the special case of $\alpha(x, t)=0$, the ELDG method becomes the classical RKDG method \cite{cockburn2001runge}. In the general setting that $\alpha$ approximates (but not exactly equals) $a$, the ELDG method enables larger time step constraint for stability than the classical DG scheme. One can compare the time step constraint \eqref{eq: time_constraint3} to that of a classical Eulerian DG method.
\end{remark}

\begin{remark} (Comparison to the ALE DG method)
It is interesting to note that when we put the Eulerian cells $I_j$ at $t^{n+1}$ and the upstream cells $I_j^\star$ at $t^n$ in a moving mesh setting, the formulation of ELDG  \eqref{mol_4} is the same as the ALE DG method \cite{klingenberg2017arbitrary} and the quasi-Lagrangian moving mesh discontinuous Galerkin method \cite{luo2019quasi}.
A fundamental difference between the ELDG and ALE DG methods is that the latter one is formulated based on a set of moving mesh, whereas the ELDG method in this paper is based on a fixed set of mesh. As a result, the ELDG method avoid the complication of mesh distortion as in an ALE DG method. In fact, the ELDG method can be viewed as a combination of SLDG algorithm in evaluating $\tilde{\mathbf{U} }^n $ and an ALE DG method in updating solutions from $\tilde{ \mathbf{U} }^n $ to $\tilde{ \mathbf{U}  }^{n+1} $.
\end{remark}
\begin{remark} (Empirical time step constraint for stability)
\label{rem:2.5}
Observe that the proposed ELDG formulation has a similar spirit to applying the RKDG method \cite{cockburn1989tvb} to 1D problems with the flux term $F = (a-\alpha)u$, thus an empirical time step stability constraint of the proposed ELDG method is
\begin{equation}
\label{eq: time_constraint2}
\Delta t \le \frac{\Delta x}{(2k+1)\max |a(x, t)-\alpha(x, t)|},
\end{equation}
with $k$ being the polynomial degree of the DG method. Combine this with \eqref{eq: time_constraint} gives
\begin{equation}
\label{eq: time_constraint3}
\Delta t \le  \frac{\Delta x}{\max\{(2k+1)\max |a(x, t)-\alpha(x, t)|,  a(x_{j+\frac12}, t^{n+1})- a(x_{j-\frac12}, t^{n+1})  \}}.
\end{equation}
For a smooth function $a$, from the construction of $\alpha$ function as previously described and by Taylor expansions, we have $\alpha - a = \mathcal{O}(\Delta t) + \mathcal{O}(\Delta x^2)$. Combining this estimate with \eqref{eq: time_constraint3} give the time step constraint for stability of ELDG
$$\Delta t \sim \Delta x^{\frac12}.$$
This is consistent with our numerical observations presented in Section~\ref{section:numerical}.
\end{remark}
\begin{remark} (Stability analysis in a simplified setting)
Stability analysis and error estimates of the proposed ELDG method solving a simplified linear equation $u_t + u_x = 0$ with $\alpha(x, t)$ for the adjoint problem being a constant $\alpha \neq 1$ could be obtained by the stability of an $L^2$ projection as in an SLDG scheme \cite{qiu2011positivity}, together with the stability of a fully discrete ALE DG method \cite{zhou2019stability}. A rigorous analysis is subject to further investigation.
\end{remark}

\begin{remark}
In our algorithm description above, $\alpha (x, t^{n+1})$ is constructed as a linear function interpolating $a(x, t)$ at cell boundaries. Alternatively, for \eqref{scalar1d}, one can track characteristics from cell boundaries at $t^{n+1}$, i.e. from $(x_{j\pm1/2}, t^{n+1})$ find their characteristics feet $(x^\star_{j\pm1/2}, t^{n+1})$. Then $\alpha(x_{j\pm1/2}, t^{n+1})$ can be obtained as the slope of the straight time connecting $(x_{j\pm1/2}, t^{n+1})$ and  $(x^\star_{j\pm1/2}, t^{n+1})$, i.e. $\alpha(x_{j\pm1/2}, t^{n+1}) = \frac{x_{j\pm1/2}-x^\star_{j\pm1/2}}{\Delta t}$.
We name the ELDG scheme with such construction of $\alpha$ function as `ELDG-ST2', and the ELDG scheme with $\alpha(x,t)$ defined by  eq.~\eqref{eq: alpha_def1} and \eqref{eq: alpha_def2}  as `ELDG-ST1' in later parts of this paper.
\end{remark}


\section{The ELDG algorithm for 2D transport problems.} \label{section:2d}

The design of the 2D ELDG algorithm shares a similar spirit as the 1D case. We consider a linear transport equation
\begin{equation}
 u_t  +  (a(x, y, t)u)_x  + (b(x, y, t)u)_y = 0, (x,y)\in \Omega.
 \label{2d_hcl}
\end{equation}
For simplicity, we assume the computational domain $\Omega$ is rectangular, boundary conditions are periodic, and   the velocity field $(a(x, y, t),b(x, y, t))$ is a continuous function of space and time.
We partition the domain $\Omega$  by a set of non-overlapping rectangular elements $A_j, j = 1,\cdots, J$,
and define the finite dimensional DG approximation space,
$V_h^k = \{ v_h:  v_h|_{A_j} \in P^k(A_j) \}$, where $P^k(A_j)$ denotes the set of polynomials of degree at most $k$
over $A_j = [x_j^l,x_j^r]\times[ y_j^b , y_j^t ]$ with element center $\left(x_j= \frac{x_j^l+x_j^r}{2} , y_{j} = \frac{ y_j^b+y_j^t }{2} \right)$ and sizes, $\Delta x_j=x_j^r - x_j^l$, $\Delta y_j = y_j^t - y_j^b$.
Let $n_k$ be the dimension of $P^k(A_j)$.

\smallskip
\noindent
{\bf (1) A modified adjoint problem for  the 2D transport problem.}
To derive a 2D  ELDG formulation, we consider a modified adjoint problem at $\tilde{A}_j(t)$ on the time interval $t\in[t^n, t^{n+1}]$:
\begin{equation}
\psi_t + \alpha(x,y,t) \psi_x + \beta(x,y,t) \psi_y =0, \quad \psi(x, y, t=t^{n+1}) = \Psi(x, y) \in P^k(A_j),
\label{eq: 2d_adjoint}
\end{equation}
where $(\alpha, \beta)$ are bilinear functions on $A_j$ at $t^{n+1}$ defined as described below. Notation-wise, we let $\tilde{A}_j(t), t \in [t^n, t^{n+1}]$  be the dynamic characteristic element of the modified adjoint problem \eqref{eq: 2d_adjoint}  with $(\tilde{x}(t), \tilde{y}(t))\in \tilde{A}_j(t)$ that satisfies \eqref{eq: char_adjoint2} emanating from $(x,y)$ of $A_j$ at $t^{n+1}$. We also let ${A}^\star_j \doteq \tilde{A}_j(t^n)$ be the upstream cell of $A_j$ at $t^n$ and let $\Omega_j$ be the region of which $(x, y, t) \in \tilde{A}_j(t) \times [t^n, t^{n+1}]$.
\begin{enumerate}
\item {\bf On $A_j$ at $t^{n+1}$.}
Let $\alpha(x, y, t^{n+1})$ and $\beta(x, y, t^{n+1}) \in Q^1(x, y)$ interpolate $a$ and $b$ functions respectively at four vertices of $A_j$,
e.g.
\begin{align}
\label{eq: a2d}
\alpha(x_{j}^l, y_{j}^b, t^{n+1}) =a(x_{j}^l, y_{j}^b, t^{n+1}), \quad
\alpha(x_{j}^l, y_{j}^t, t^{n+1}) =a(x_{j}^l, y_{j}^t, t^{n+1}),  \\
\alpha(x_{j}^r, y_{j}^b, t^{n+1}) =a(x_{j}^r, y_{j}^b, t^{n+1}), \quad
\alpha(x_{j}^r, y_{j}^t, t^{n+1}) =a(x_{j}^r, y_{j}^t, t^{n+1}). \nonumber
\end{align}
Similarly, $\beta$ is a bilinear function interpolating $b$ at four vertices $(x_{j}^{l}, y_{j}^{b})$, $(x_{j}^{l}, y_{j}^{t})$, $(x_{j}^{r}, y_{j}^{b})$, $(x_{j}^{r}, y_{j}^{t})$.
\item {\bf On $\tilde{A}_j(t)$ at $t\in[t^n, t^{n+1})$.}  Along characteristic lines of the adjoint problem \eqref{eq: 2d_adjoint} emanating from  any point $(\xi, \eta)\in A_j$ at $t^{n+1}$,  with
 $$\tilde{x}(t; (\xi, \eta, t^{n+1})), \tilde{y}(t; (\xi, \eta, t^{n+1}))$$
 satisfy the following equations,
\begin{equation}
\label{eq: char_adjoint2}
 \frac{d}{d t}\tilde{x}(t; (\xi, \eta, t^{n+1})) = \alpha(\xi, \eta, t^{n+1}), \quad
\frac{d}{d t}\tilde{y}(t;  (\xi, \eta, t^{n+1})) = \beta(\xi, \eta, t^{n+1}),
\end{equation}
from which one have
\begin{equation}
\tilde{x}(\tau; (\xi, \eta, t^{n+1})) = \xi - \alpha(\xi, \eta, t^{n+1}) (t^{n+1}-\tau) \in Q^1(\xi, \eta),
\label{eq: tildex_2d}
\end{equation}
\begin{equation}
\tilde{y}(\tau; (\xi, \eta, t^{n+1})) = \eta - \beta(\xi, \eta, t^{n+1}) (t^{n+1}-\tau) \in Q^1(\xi, \eta),
\label{eq: tildey_2d}
\end{equation}
with the Jacobian
\begin{equation}
J(\xi, \eta, \tau) = \frac{\partial (\tilde{x}, \tilde{y})}{\partial(\xi, \eta)}(\tau) =
   \left(
   \begin{array}{cc}
    1-\frac{\partial \alpha}{\partial \xi}(t^{n+1}-\tau) & \frac{\partial \alpha}{\partial \eta}(t^{n+1}-\tau)\\
    -\frac{\partial \beta}{\partial \xi}(t^{n+1}-\tau) & 1-\frac{\partial \beta}{\partial \eta}(t^{n+1}-\tau)
   \end{array}
   \right).
\label{eq: 2d_jacobian}
\end{equation}
Then we let, for $t\in[t^n, t^{n+1}]$, and $(\tilde{x}, \tilde{y})\in \tilde{A}_j(t)$,
\begin{equation}
\alpha(\tilde{x}(t; (\xi, \eta, t^{n+1})), \tilde{y}(t; (\xi, \eta, t^{n+1})), t)=\alpha(\xi, \eta, t^{n+1}),
\label{eq: alpha}
\end{equation}
\begin{equation}
\beta(\tilde{x}(t; (\xi, \eta, t^{n+1})), \tilde{y}(t; (\xi, \eta, t^{n+1})), t)=\beta(\xi, \eta, t^{n+1}).
\label{eq: beta}
\end{equation}
It can be easily checked that, \eqref{eq: char_adjoint2} are the characteristics equations for the modified adjoint problem \eqref{eq: 2d_adjoint} with $\alpha$ and $\beta$ functions defined by eq.~\eqref{eq: alpha} and \eqref{eq: beta}.
For the modified adjoint problem, the solution $\psi$ stays constant along characteristics, therefore we have
\begin{equation}
\psi(\tilde{x}(\tau;  (\xi, \eta, t^{n+1})),  \tilde{y}(\tau;  (\xi, \eta, t^{n+1})),\tau) = \Psi(\xi, \eta) \in P^k(A_j), \quad \forall \tau\in[t^n, t^{n+1}].
\label{testfunction2d}
\end{equation}
\end{enumerate}

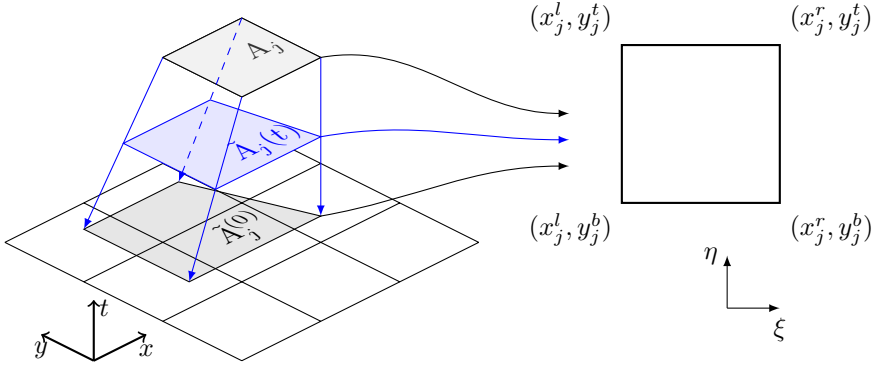
\begin{figure}
\centering
\begin{tikzpicture}[scale=0.7]
       \begin{scope}[
            yshift= 0,every node/.append style={
            yslant=0.5,xslant=-1},yslant=0.5,xslant=-1
            ]
 \draw[black,fill=black!10] (1,4)  -- ( 1,2)  -- (3.5,2) -- (2.8 ,4)       -- cycle;
  \node[black!100, rotate=0] (a) at ( 2.5,2.5) {  $ \tilde{A}_j^{(0)} $ };
    \coordinate (A) at (1,2);
     \coordinate (B) at (3.5,2);
      \coordinate (C) at (2.8,4);
       \coordinate (D) at (1,4);
    \end{scope}
    \begin{scope}[
            yshift=0,every node/.append style={
            yslant=0.5,xslant=-1},yslant=0.5,xslant=-1
            ]
          \draw[step=1.5cm, black] (0,0) grid (4.5,4.5); 

        \coordinate (V) at (3.5,  2);

    \end{scope}

        \begin{scope}[
            yshift= 50,every node/.append style={
            yslant=0.5,xslant=-1},yslant=0.5,xslant=-1
            ]

  \draw[blue,fill=blue!10] (0.5*1+0.5*1.5,  0.5*2+0.5*1.5) -- ( 0.5*3.5+0.5*3, 0.5*2+0.5*1.5)  -- ( 2.8*0.5+0.5*3,4*0.5+3*0.5) -- (1.5*0.5+1*0.5 ,3*0.5+4*0.5)       -- cycle;
  \node[blue!100, rotate=0] (a) at ( 2.6,2.2) {  $\tilde{A}_j(t)$ };

    \coordinate (U) at (0.5*3.5+0.5*3, 0.5*2+0.5*1.5);

    \end{scope}
        \begin{scope}[
            yshift= 100,every node/.append style={
            yslant=0.5,xslant=-1},yslant=0.5,xslant=-1
            ]
          \draw[step=1.5cm, black ] (1.5,1.5) grid (3,3); 

  \draw[black,fill=black!5] (1.5,1.5)  -- ( 3,1.5)  -- (3,3) -- (1.5 ,3)       -- cycle;
  \node[black!80, rotate=0] (a) at ( 2.6,2.2) {  $ A_j $ };

    \draw[-latex,blue]  (1.5,1.5) node[right,scale=1.0]{}--(A)  node[] {};
    \draw[-latex, blue]  (3,1.5) node[right,scale=1.0]{}--(B)  node[] {};
    \draw[-latex,dashed,blue]  (3,3) node[right,scale=1.0]{}--(C)  node[] {};
    \draw[-latex, blue]  (1.5,3) node[right,scale=1.0]{}--(D)  node[] {};

\coordinate (W) at (3,  1.5);

    \end{scope}
    \begin{scope}[
            yshift= 0,xshift= -80,every node/.append style={
            yslant=0.5,xslant=-1},yslant=0.5,xslant=-1
            ]
     \coordinate (O) at (0,0,0);
  \draw[thick,->] (0,0,0) -- (1,0,0) node[anchor=north east]{ };
  \draw[thick,->] (0,0,0) -- (0,1,0) node[anchor=north west]{ };
  \draw[thick,->] (0,0,0) -- (0,0,-3) node[anchor=north]{ };
  \end{scope}

      \begin{scope}[
            yshift= 0,xshift= -80
            ]
            \node[ rotate=0] (a) at ( -1,0.2) {  $y$ };
            \node[ rotate=0] (a) at ( 1 ,0.2) {  $x$ };
            \node[rotate=0] (a) at ( 0.2 ,1) {  $t$ };

  \end{scope}

  \begin{scope}
  [
            yshift= 0,xshift= 120
            ]

    \draw[black,thick ] (3,3)node[below left,scale=1.0]{  $(x_j^l,y_j^b )$}
     -- ( 6,3) node[below right,scale=1.0]{   $(x_j^r,y_j^b)$}
      -- (6,6)node[above right,scale=1.0]{   $(x_{j}^r,y_{j}^t)$}
       -- ( 3 ,6)  node[above left,scale=1.0]{  $(x_{j}^l,y_{j}^t )$}     -- cycle;

\draw[-latex] (5,1) -- (5,1+1) node[left]{ $\eta$ };
  \draw[-latex] (5,1) -- (5+1,1) node[below]{ $\xi$ };

\draw[-latex, blue]  (U)  node[] {} to[out=10,in=180] (2,4.2) node[right,scale=1.0]{} ;
   \draw[-latex, black]  (W)  node[] {} to[out=10,in=180] (2,4.7) node[right,scale=1.0]{} ;
\draw[-latex, black]  (V)  node[] {} to[out=10,in=180] (2,3.7) node[right,scale=1.0]{} ;
  \end{scope}

\end{tikzpicture}
\caption{Illustration for the mapping between dynamic element $\tilde{A}_j(t)$ (left) and the iso-parametric element (right).   }
\label{2dmap}
\end{figure}

Next we introduce a few notations and useful equalities \cite{ciarlet1988mathematical,persson2009discontinuous} regarding the coordinate transformation defined by \eqref{eq: tildex_2d}-\eqref{eq: tildey_2d} .
 \begin{equation}
    dxdy = \det( J(\xi,\eta) )  d\xi d\eta,
 \end{equation}
\begin{equation}
\nabla_{x,y} \psi(x,y) =  J(\xi,\eta) ^{-T} \nabla_{\xi,\eta} \Psi(\xi,\eta),
\end{equation}
\begin{equation}
\mathbf{n} dS = \det( J(\xi,\eta) ) J(\xi,\eta)^{-T} \breve{ \mathbf{n}  } d\breve{S},
\end{equation}
where $dS$ and $d\breve{S}$ are the infinitesimal boundaries of the dynamic element and the isoparametric element, respectively and their corresponding normal vectors are $\mathbf{n}$ and $\breve{ \mathbf{n}  }$. The inverse of the Jacobian is given by
\begin{equation}
 J( \xi , \eta )^{-1}
 =
 \frac{1}{ |\det( J(\xi,\eta) )| }
 \left(
\begin{array}{cc}
 \tilde{y}_{\eta}  &  - \tilde{x}_{\eta} \\
 -\tilde{y}_{\xi}  &  \tilde{x}_{\xi}
\end{array}
 \right).
\end{equation}
We assume the determinant of the Jacobian $\det( J(\xi,\eta) )$ is positive; if the Jacobian is negative, it indicates the distortion of upstream cells. In such a situation, the time stepping size should be reduced by using the adaptive time stepping algorithm \cite{cai2019exp}.

\smallskip
\noindent
{\bf (2) Semi-discrete ELDG scheme formulation.}
Integrating $\eqref{2d_hcl}\cdot \psi+\eqref{eq: 2d_adjoint} \cdot u$ over $\Omega_j$, we have
\begin{equation}
  \int_{ \Omega_j }  \left[ \eqref{2d_hcl}\cdot \psi+\eqref{eq: 2d_adjoint} \cdot u \right]   dxdydt =0.
\end{equation}
Then,
\begin{align}
0 =&\int_{t^n }^{ t^{n+1} } \int_{ \tilde{A}_j(t)   } ( u_t\psi+u\psi_t ) dxdydt  \nonumber \\
   & +  \int_{t^n }^{ t^{n+1} } \int_{ \tilde{A}_j(t)   } (  (a u)_x\psi + \alpha\psi_x u + (b u)_y\psi+\beta \psi_y u  )  dxdy dt \nonumber\\
  =&\int_{t^n }^{ t^{n+1} } \left[ \int_{ \tilde{A}_j(t)   } ( u\psi )_t dxdydt  +    \int_{ \tilde{A}_j(t)   } (  (a u)_x\psi + \alpha\psi_x u + (b u)_y\psi+\beta \psi_y u  )  dxdy dt \right] \nonumber \\
  =&\int_{ t^n }^{ t^{n+1} } \left[ \frac{d}{dt} \int_{ \tilde{A}_j(t) } u\psi dxdy -\int_{ \partial \tilde{A}_j(t) } u\psi
   \left(
   \begin{array}{c}
     \alpha \\ \beta
   \end{array}
   \right)
   \cdot\mathbf{n} dS \right] dt \nonumber \\
   & + \int_{ t^n }^{ t^{n+1} } \left[ \int_{ \tilde{A}_j(t) } \nabla\cdot
   \left(
   \begin{array}{c}
    au  \\ bu
   \end{array}
   \right) \psi dxdy
   +
   \int_{ \tilde{A}_j(t) }
   \left(
   \begin{array}{c}
    \alpha\\  \beta
   \end{array}
   \right)
   \cdot
   \nabla \psi
   u dxdy \right] dt  \nonumber\\
    =&\int_{ t^n }^{ t^{n+1} } \left[ \frac{d}{dt} \int_{ \tilde{A}_j(t) } u\psi dxdy +\int_{ \partial \tilde{A}_j(t) } \psi
 \mathbf{F}
    \cdot\mathbf{n} dS
- \int_{\tilde{A}_j(t) }
  \mathbf{F} \cdot
   \nabla \psi dxdy
    \right] dt, \nonumber \\
\label{semi_SL0}
\end{align}
with \begin{equation}
\mathbf{F}(u,x,y,t)
=
\left(
\begin{array}{c}
(a(x, y, t) - \alpha(x,y,t))u \\
(b(x, y, t) - \beta(x,y,t))u
\end{array}
\right),
\end{equation}
in which the Leibniz-Reynolds transport theorem and the divergence Theorem are used for the above derivation.
The time differential version of eq.~\eqref{semi_SL0} can be written as
\begin{align}
 &  \frac{d}{dt} \int_{ \tilde{A}_j(t) } u\psi dxdy     =  -  \int_{ \partial \tilde{A}_j(t) } \psi
    \mathbf{F} \cdot \mathbf{n}dS +
    \int_{ \tilde{A}_j(t)  }
  \mathbf{F}
   \cdot
   \nabla \psi dxdy.
   \label{semi_SL1}
\end{align}

As the 1D case, we map the coordinate of $(x, y)\in \tilde{A}_j(t)$ to a reference cell of $(\xi, \eta)\in A_j$ as shown in Figure \ref{2dmap}.
Then we rewrite eq.~\eqref{semi_SL1}   as
\begin{align}
 &  \frac{d}{dt} \int_{ {A}_j } u(\tilde{x}(t, (\xi, \eta, t^{n+1})), \tilde{y}(t, (\xi, \eta, t^{n+1})), t) \Psi(\xi, \eta) \det(J(\xi, \eta, t)) d\xi d\eta     \nonumber \\
 & =-    \int_{ \partial A_j } \Psi(\xi,\eta) \mathbf{F}\cdot \left(  \det( J(\xi,\eta, t) ) J(\xi,\eta, t)^{-T} \breve{ \mathbf{n}  } \right) d\breve{S} \nonumber \\
  &  + \int_{{A}_j}
  \mathbf{F}
   \cdot
   (J(\xi, \eta, t)^{-T} \nabla_{\xi, \eta} \Psi) \det(J(\xi, \eta, t)) d\xi d\eta
      .
       \label{eq: semi_SL3}
\end{align}
Notice that in equation~\eqref{eq: semi_SL3}, functions are all in the $(\xi, \eta)$ coordinate, and can be evolved by the method-of-lines approach, e.g. using explicit SSP RK methods. $\Psi(\xi, \eta)$ function stays as the same polynomial in the $(\xi, \eta)$ coordinate for all $t\in[t^n, t^{n+1}]$ by the design of our adjoint problem, see eq.~\eqref{testfunction2d}.

%
We let the approximate solution of $u(\tilde{x}(t, (\xi, \eta, t^{n+1})), \tilde{y}(t, (\xi, \eta, t^{n+1})), t)$ be written in the $(\xi, \eta)$ coordinate as follows,
 \begin{equation}
  u_h(\xi,\eta,t) = \sum_{p=1}^{ n_k } \breve{u}_p(t) \Psi_p(\xi,\eta),
  \label{eq: u_h_repres}
 \end{equation}
 where bases $\Psi_p(\xi,\eta),p=1,\cdots, n_k$ expands the space of $P^k(A_j)$, for implementation.
 For the ELDG scheme, we look for $u_h$ in the above form satisfying
 \begin{align}
 && \frac{d}{dt} \int_{ {A}_j } u_h  \Psi_p  \det(J(\xi, \eta, t)) d\xi d\eta
+    \int_{ \partial A_j } \Psi_p \hat{\mathbf{F}}\cdot \left(  \det( J(\xi,\eta, t) ) J(\xi,\eta, t)^{-T} \breve{ \mathbf{n}  } \right) d\breve{S}
 \nonumber\\ &&   - \int_{{A}_j}
  \mathbf{F}
   \cdot
   (J(\xi, \eta, t)^{-T} \nabla_{\xi, \eta} \Psi) \det(J(\xi, \eta, t)) d\xi d\eta = 0.
   \label{eq: 2d_ELDG_semi}
 \end{align}
Here $\hat{\mathbf{F}}$ in the second term is a monotone numerical flux, an example of which is the Lax-Friedrich flux, and the line and volume integral in the second and third terms could be performed by proper high order quadrature rules as in a standard RK DG scheme.
Then the coefficients $\textbf{u} = ( \breve{u}_1,\breve{u}_2, \cdots,\breve{u}_{n_k} )^T$ in \eqref{eq: u_h_repres} satisfies a system of ODEs,
 \begin{equation}
 \label{eq: semi-2d}
 \frac{d}{dt} \left(\textbf{M}(t)\textbf{u}(t)\right)=L(\textbf{u}(t)),
 \end{equation}
where  the mass matrix $\textbf{M}$ is of size $n_k$ by $n_k$ and its entries are
  \begin{equation*}
    M_{pq}(t)
    = \int_{A_j}
                 \Psi_p(\xi,\eta) \Psi_q(\xi,\eta) \det( J(\xi,\eta, t) ) d\xi d\eta,
  \end{equation*}
 and $L(\textbf{u}(t))$ is the RHS vector from the evaluation of the other terms in \eqref{eq: 2d_ELDG_semi}.



\noindent{\bf (3) RK time discretization and fully discrete scheme.}
The semi-discrete scheme \eqref{eq: semi-2d} can be discretized by applying an explicit RK time discretization with the initial condition
\begin{equation}
\textbf{M}(t^n) \textbf{u}^n = \int_{{A}_j^\star } u^n_h(x, y) \psi(x, y, t^n) dxdy,
\label{semi_2d_init}
\end{equation}
being evaluated by a 2D SLDG procedure \cite{cai2016high}.
Below we provide a flow chart of the fully discrete 2D algorithm described above.
\begin{description}
  \item[Step 1.] Construct $\alpha(x,y,t)$ and $\beta(x,y,t)$  for $(x, y, t)\in \tilde{A}_j(t)\times[t^n, t^{n+1}]$ by first constructing
  $$\alpha(x, y, t^{n+1}), \beta(x, y, t^{n+1})\in Q^1(x, y), \quad (x, y)\in A_j,$$
 interpolating $a(x, y, t^{n+1}), b(x, y, t^{n+1})$ respectively at four vertices of $A_j$; then these $\alpha$ and $\beta$ functions are constructed by following \eqref{eq: alpha}-\eqref{eq: beta} for $t\in[t^n, t^{n+1})$.  In particular, one first find $(\xi, \eta)$ for $(\tilde{x}$, $\tilde{y})$ from \eqref{eq: tildex_2d}-\eqref{eq: tildey_2d}; then the $\alpha(\tilde{x}$, $\tilde{y}, t)$ and $\beta(\tilde{x}$, $\tilde{y}, t)$ are defined following \eqref{eq: alpha}-\eqref{eq: beta}. Note that, while $(\tilde{x}, \tilde{y})$ is a bilinear function of $(\xi, \eta)$, the same statement does not hold for the inverse mapping. Figure~\ref{2dmap} illustrates 2D transformation between $(\xi, \eta)\in A_j$ and $(x, y)\in \tilde{A}_j(t)$ for some $t\in[t^n, t^{n+1}]$.

 \item[Step 2.]

 Set up dynamic elements $\tilde{A}_j^{(l)}, l=0,\cdots,s,$ for each immediate stage of the RK method, and compute the corresponding Jacobian of the transformation  $J=\frac{\partial(x, y)}{\partial(\xi, \eta)}$, $ J(\xi, \eta, \tau)^{-1}$ in  \eqref{eq: semi_SL3}; these quantities can be precomputed as functions of $(\xi, \eta, t^{(l)})$.

\item[Step 3.] Perform the SLDG algorithm in \cite{cai2016high} to get  the initial condition of \eqref{semi_2d_init}.
   Notice that
   since the mapping $(x(\xi,\eta), y(\xi,\eta) )$ in \eqref{eq: tildex_2d}-\eqref{eq: tildey_2d} is not affine, it is not as straighforward to find the inverse mapping of $(\xi(x,y), \eta(x,y) )$ as the 1D problem. Some approximation, as is done in \cite{cai2016high}, has to be performed in order to obtain $\psi(x, y, t^n)$.

   \item[Step 4.]  An SSP RK method is applied to \eqref{eq: semi-2d}.
     In particular, at the $l^{th}$ RK stage, $\textbf{M}^{(l)} \textbf{u}^{(l)}$ is first being updated, then $\textbf{u}^{(l)}$ is computed by applying $(\textbf{M}^{(l)})^{-1}$; finally $\textbf{u}^{(l)}$ as the degree of freedom in $(\xi, \eta)$ coordinate are being used to evaluate the RHS of \eqref{eq: semi-2d} for future RK stages.
\end{description}
\begin{remark} (Quadrilateral shape of upstream cells)
The fact that $\alpha(x, y, t^{n+1})$ and $\beta(x, y, t^{n+1})$ functions are in $Q^1(A_j)$ in the modified adjoint problem ensures the quadrilateral shape of upstream cells. This avoids the need to use quadratic curves to approximate upstream cells in achieving high order spatial accuracy in the original SLDG algorithm \cite{cai2016high}. An example of such is the swirling deformation example as shown in the numerical section.
\end{remark}
\begin{remark} (Assumption on the velocity field)
For the scope and applications of our current paper, we work with the velocity fields $(a(x, y, t), b(x, y, t))$ that are smooth enough and divergence free. The proposed ELDG formulation works for general non-divergence free velocity field as long as the Jacobian of the transformation is always positive.
\end{remark}

\section{ELDG method with the exponential integrators for nonlinear Vlasov dynamics} \label{section:nonlinear}

The proposed ELDG method for linear transport problems can be applied to solve nonlinear models such as Vlasov models, via combining with the Runge-Kutta exponential integrator method in \cite{celledoni2003commutator, cai2019exp}. We will denote such a method as ELDG-RKEI. Below we first present the nonlinear Vlasov-Poisson, the guiding center Vlasov models
as well as    the 2D incompressible Euler equations; and then present a second order and a third order ELDG-RKEI method.

\noindent
{\bf The nonlinear Vlasov-Poisson system} reads as follows,
  \begin{equation}
f_t + vf_x + E(x,t) f_v =0,
\label{vlasov}
\end{equation}
\begin{equation}
E(x,t) = -\phi_x, \
-\phi_{xx}(x,t) = \rho(x,t),
\label{poisson}
\end{equation}
where the electron distribution function $f(x,v,t)$ is the probability distribution function in the phase space $(x,v)\in\Omega_x \times \mathbb{R}$ describing the probability of finding a particle with velocity $v$ at position $x$ and at time $t$.
The electric field $E=-\phi_x$, where the self-consistent electrostatic potential $\phi$ is determined by the Poisson's equation \eqref{poisson}. $\rho(x,t) = \int_{ \mathbb{R} } f(x,v,t) dv -1$ denotes charge density, with the assumption that infinitely massive ions are uniformly distributed in the background.

\noindent
{\bf The guiding center Vlasov model} describes a highly magnetized plasma in the transverse plane of a tokamak \cite{shoucri1981two,crouseilles2009conservative}, and reads as follows:
\begin{equation}
\rho_t + \nabla \cdot(  \mathbf{E}^{\bot} \rho ) = 0,
\label{guiding}
\end{equation}
\begin{equation}
-\Delta \Phi = \rho,\
\mathbf{E}^{\bot} =  ( -\Phi_y,\Phi_x ),
\label{poisson2d}
\end{equation}
where the unknown variable $\rho$ denotes the charge density of the plasma, and the electric field $\mathbf{E}$ depends on $\rho$ via the Poisson equation.

\noindent
{\bf  The  2D incompressible Euler in the vorticity-stream function} reads as follows,
\begin{equation}
\omega_t + \nabla\cdot( \mathbf{u} \omega  ) =0,
\label{Euler}
\end{equation}
\begin{equation}
 \Delta \Phi = \omega, \
 \mathbf{u} = -( -\Phi_y , \Phi_x ),
\end{equation}
where $\mathbf{u}$ is the velocity field, $\omega$ is the vorticity of the fluid, and $\Phi$ is the stream-function determined by Poisson's equation.

The above three models can be written in the form of \eqref{general_nonlinear}.
In \cite{celledoni2003commutator, celledoni2009semi, cai2019exp}, the exponential integrator method is applied to solve nonlinear time-dependent problems \eqref{general_nonlinear}, by decomposing the nonlinear dynamics into the composition of a sequence of linearized transport problems to achieve high order temporal accuracy.
We denote the ELDG procedure of updating the solution of linearized equation from $t^*$ to $t^*+\Delta t$ with frozen velocity field $\mathbf{P}( u^*;\mathbf{x},t^* )$
\begin{equation}
\begin{cases}
u_t + \nabla \cdot( \mathbf{P}( u^*;\mathbf{x},t^* ) u ) = 0,  \\
u(t^*) = u^*,
\end{cases}
\end{equation}
as
\begin{equation}
ELDG(\mathbf{P}(u^{*}; \mathbf{x},t^* ),\Delta t )(u^*).
\end{equation}
When a second order RKEI scheme is used with the ELDG update of linearized solution, one has
\begin{align*}
u^{(1)}  &  =  u^n \nonumber \\
u^{(2)}  &  =  ELDG \left( \frac12 \mathbf{P}( u^{(1)} ),  \Delta t    \right)  u^{(1)} \nonumber \\
u^{n+1}  &  =  ELDG \left(  \mathbf{P}( u^{(2)} )     ,
\Delta t \right)  u^{(1)}.
\end{align*}
We name such scheme `ELDG-CF2' \cite{cai2019exp}, in which `CF2' refers to the above second order RKEI scheme.
When a third order RKEI scheme is used with the ELDG update of linearized solution, one has
\begin{align*}
u^{(1)}  &  =  u^n \nonumber \\
u^{(2)}  &  =  ELDG \left( \frac13 \mathbf{P}( u^{(1)} ),  \Delta t    \right)  u^{(1)} \nonumber \\
u^{(3)}  &  =  ELDG \left(   \frac23 \mathbf{P}( u^{(2)} ) , \Delta t    \right) u^{(1)} \nonumber \\
u^{n+1}  &  =  ELDG \left( - \frac{1}{12} \mathbf{P}( u^{(1)} )      + \frac{3}{4} \mathbf{P}( u^{(3)} ),
\Delta t \right)  u^{(2)}.
\end{align*}
We name such scheme `ELDG-CF3C03' \cite{cai2019exp}, in which `CF3C03' refers to the above third order RKEI scheme.
We refer to \cite{cai2019exp} for more details regarding implementation.
In the nonlinear Vlasov models LDG schemes \cite{arnold2002unified, cockburn1998local, castillo2000priori,James} are adopted to solve the elliptic field equations \eqref{poisson} and \eqref{poisson2d}.

\section{Numerical results}
\label{section:numerical}

In this section, we perform numerical experiments for linear transport problems as well as the nonlinear Vlasov models.
To showcase the proposed method, we perform the following studies:
(1) the convergence of spatial discretization by using small enough time stepping size;
(2) we vary $CFL$ to study the temporal convergence and numerical stability with a well resolved spatial mesh;
(3) we present snapshots of numerical solutions in a long time;
(4) we numerically track the time history of invariants, such as mass and energy.

The ELDG method presented below is the ELDG-ST1 method, unless otherwise noted. When needed, we use the $k+1$-th order RK for tracing characteristic lines.
We set the time step for 1D and 2D problems as
 \begin{equation}
 \Delta t= CFL\Delta x   \ \text{and} \  \Delta t=\frac{CFL}{ \frac{a}{\Delta x}+\frac{b}{\Delta y} },
 \end{equation}
 respectively; here $a$ and $b$ are maximum transport speeds in $x$ and $y$ directions, respectively.
For some test cases, we also present the SLDG \cite{cai2016high,cai2019exp} and classical RKDG methods for comparison purpose.

\subsection{1D linear transport problems}
\begin{example}
(1D linear transport equation with constant coefficient.)
We start with the following 1D transport equation
\begin{equation}
u_t + u_x =0, \ x\in[0,2\pi],
\end{equation}
with the smooth initial data $u(x,0)=\sin(x)$ and exact solution $u(x,t)=\sin(x-t)$.
For the constant coefficient problem, the proposed ELDG method, if using the exact velocity field, is the same as SLDG.
Here we perturb the velocity at cell boundaries for the modified adjoint problem to be $\alpha(x_{j+\frac12}) = 1+\sin(x_{j+\frac12})\Delta x$.

Table \ref{linear1d_spatial} reports the spatial accuracies of the ELDG, SLDG and RKDG methods for this example with the same time stepping size. The proposed ELDG method is found to be as accurate as the SLDG and RKDG methods.
We vary time stepping size, with fixed well-resolved spatial meshes, and plot error vs. $CFL$ in Figure \ref{linear_stability} for ELDG and SLDG $P^1$ (left) and $P^2$ (right) schemes at a long time $T=100$. For the ELDG scheme, the time-stepping constraint can be found to be $\Delta t \leq \frac{1}{(2k+1)\Delta x} \Delta x$ from the perturbation of velocity field and \eqref{eq: time_constraint2}; hence
$$CFL_{\text{upper}}= \frac{1}{  (2k+1)\Delta x }, $$
for $P^k$ ELDG schemes. They are shown as dashed lines in the figure. It is observed that these bounds are expected in this numerical test. The SLDG schemes are observed to be unconditionally stable. The ELDG and SLDG schemes are observed to have similar error magnitudes, when the $CFL$ is less than the stability bounds (dash lines).

\begin{table}[!ht]\footnotesize
\caption{1D linear transport equation with constant coefficient. $u_t+u_x=0$ with initial condition $u(x,0) = \sin(x)$. $T=\pi$.
We use $CFL=0.3$ and $CFL=0.18$ for all $P^1$ and $P^2$ schemes, respectively. ELDG here with the vertex perturbation. 
  }
\centering
\begin{tabular}{| c | cc  | cc| cc| }

\hline
Mesh  &{$L^1$ error} & Order    &  {$L^1$ error} & Order &  {$L^1$ error} & Order  \\
 \hline
  & \multicolumn{2}{c|}{$P^1$ RKDG} & \multicolumn{2}{c|}{$P^1$ SLDG}& \multicolumn{2}{c|}{$P^1$ ELDG}
  \\ \hline
    40 &     1.15E-03 &     -- &    6.37E-04 &    --       &     6.08E-04 &    -- \\
    80 &     2.85E-04 &     2.01 &    1.59E-04 &     2.00 &     1.55E-04 &     1.97 \\
   160 &     7.09E-05 &     2.01 &    3.90E-05 &     2.03  &     3.84E-05 &     2.02 \\
   320 &     1.77E-05 &     2.00 &     1.77E-05 &     2.00 &     9.77E-06 &     1.98 \\
\hline
  & \multicolumn{2}{c|}{$P^2$ RKDG} & \multicolumn{2}{c|}{$P^2$ SLDG}& \multicolumn{2}{c|}{$P^2$ ELDG}
  \\ \hline
    40 &    9.28E-06 & --    &    7.25E-06 &    --        &     7.69E-06 & -- \\
    80 &     1.16E-06 &     3.00 &   9.23E-07 &     2.97  &     9.45E-07 &     3.03 \\
   160 &     1.44E-07 &     3.00 &    1.17E-07 &     2.98   &    1.18E-07 &     3.00 \\
   320 &    1.80E-08 &     3.00  &     1.40E-08 &     3.06 &    1.41E-08 &     3.07 \\
\hline
\end{tabular}
\label{linear1d_spatial}
\end{table}

\begin{figure}[!ht]
\centering
\includegraphics[height=50mm]{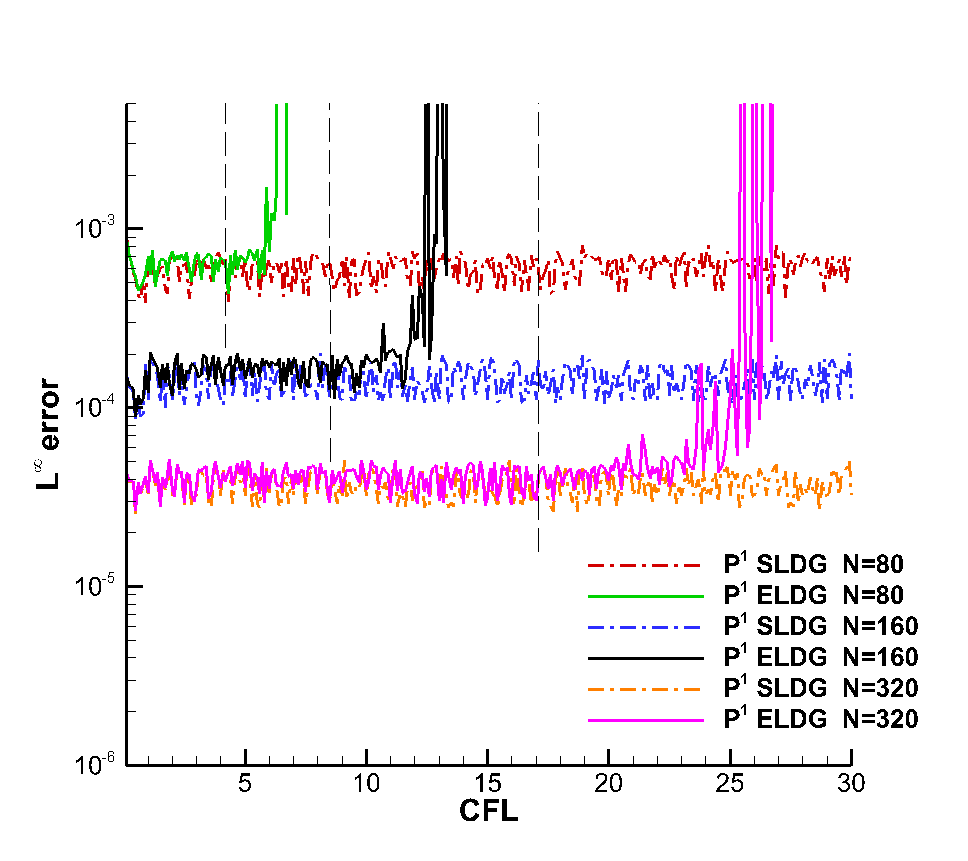}
\includegraphics[height=50mm]{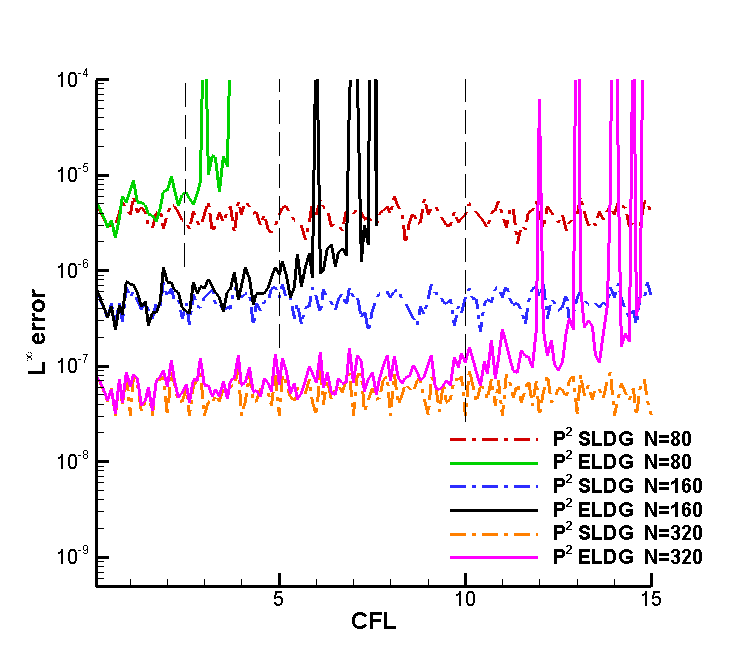}
\caption{The $L^\infty$ error versus $CFL$ of SLDG methods and ELDG methods for 1D linear transport equation with constant coefficient:  $u_t+u_x=0$ with initial condition $u(x,0) = \sin(x)$. A long time simulation is performed with $T=100$.  
The vertical long dashes from left to right are expected upper bounds of $CFL$ for stability for $P^k$ ELDG methods with meshes $80$, $160$ and $320$ respectively.
}
\label{linear_stability}
\end{figure}

\end{example}

\begin{example}
(1D transport equation with variable coefficients.)
Consider
\begin{equation}
u_t + (\sin(x)u)_x =0, \ x\in[0,2\pi]
\end{equation}
with initial condition $u(x,0)=1$ and the periodic boundary condition.
The exact solution is given by
\begin{equation}
u(x,t) =
\frac{  \sin( 2\tan^{-1}( e^{-t}\tan(\frac{x}{2}) ) )   }{  \sin(x) }.
\end{equation}

As in the previous example, the spatial convergence of RKDG, SLDG, ELDG-ST1 and ELDG-ST2 are shown in Table \ref{1dsin_example}. The expected spatial convergence orders are observed.
In Figure \ref{sin1d_time}, we plot the $L^\infty$ error versus $CFL$ of ELDG-ST1, ELDG-ST2 and SLDG schemes with $P^1$ (left) and $P^2$ (right) polynomial spaces. The following observations are made:
(1) all methods perform similarly around and before $CFL=1$, which is well above the stability constraint of the RKDG method $1/(2k+1)$; (2) after $CFL=1$ and before stability constraint of the method, the temporal convergence order is observed to be $2$ and $3$ for $P^1$ and $P^2$ respectively, corresponding to the RK method used in time integration and characteristics tracing; (3) the upper bounds of $CFL$ for stability of $P^2$ ELDG with mesh $N=80,160,320$ are around $3.5$, $5$, $7$, which increase with ratio around $\sqrt{2}$. This verifies the time step estimate $\Delta t \sim \sqrt{\Delta x}$ in Remark~\ref{rem:2.5}.

\begin{table}[!ht]
\footnotesize
\caption{1D transport equation with variable coefficients. $u_t + (\sin(x)u)_x =0$ with the initial condition $u(x,0) = 1$. $T=1$.
We use $CFL=0.3$ and $CFL=0.18$ for all $P^1$ and $P^2$ schemes, respectively.
  }
\centering
\begin{tabular}{| c | cc  | cc| cc| cc|  }
\hline
 Mesh  &{$L^1$ error} & Order    &  {$L^1$ error} & Order &  {$L^1$ error} & Order &  {$L^1$ error} & Order \\
 \hline
  & \multicolumn{2}{c|}{$P^1$ RKDG} & \multicolumn{2}{c|}{$P^1$ SLDG}& \multicolumn{2}{c|}{$P^1$ ELDG-ST1} & \multicolumn{2}{c|}{$P^1$ ELDG-ST2}
  \\ \hline
    40 &     1.30E-03 &   -- &    1.35E-03 &    --       &     1.20E-03 &    -- &     1.35E-03 &   --\\
    80 &     3.25E-04 &     2.00  &    3.56E-04 &     1.92 &     3.24E-04 &     1.89 &3.54E-04 &     1.93\\
   160 &     8.14E-05 &     2.00 &    8.95E-05 &     1.99  &     8.35E-05 &     1.96 &8.89E-05 &     1.99\\
   320 &     2.04E-05 &     2.00 &     2.31E-05 &     1.95 &     2.21E-05 &     1.92 &     2.30E-05 &     1.95\\
\hline
  & \multicolumn{2}{c|}{$P^2$ RKDG} & \multicolumn{2}{c|}{$P^2$ SLDG}& \multicolumn{2}{c|}{$P^2$ ELDG-ST1}& \multicolumn{2}{c|}{$P^2$ ELDG-ST2}
  \\ \hline
    40 &    8.11E-05 &  --       &   5.16E-05 &     --        &    6.45E-05 & -- &     5.20E-05 & --\\
    80 &    1.21E-05 &     2.74  &   6.35E-06 &     3.02  &     7.36E-06 &     3.13 &     6.36E-06 &     3.03\\
   160 &    1.79E-06 &     2.76  &   7.85E-07 &     3.02   &    8.65E-07 &     3.09&     7.87E-07 &     3.02 \\
   320 &    2.62E-07 &     2.78  &    9.61E-08 &     3.03 &    1.02E-07 &     3.08&     9.63E-08 &     3.03 \\
\hline
\end{tabular}
\label{1dsin_example}
\end{table}

\begin{figure}[!ht]
\centering
\includegraphics[height=50mm]{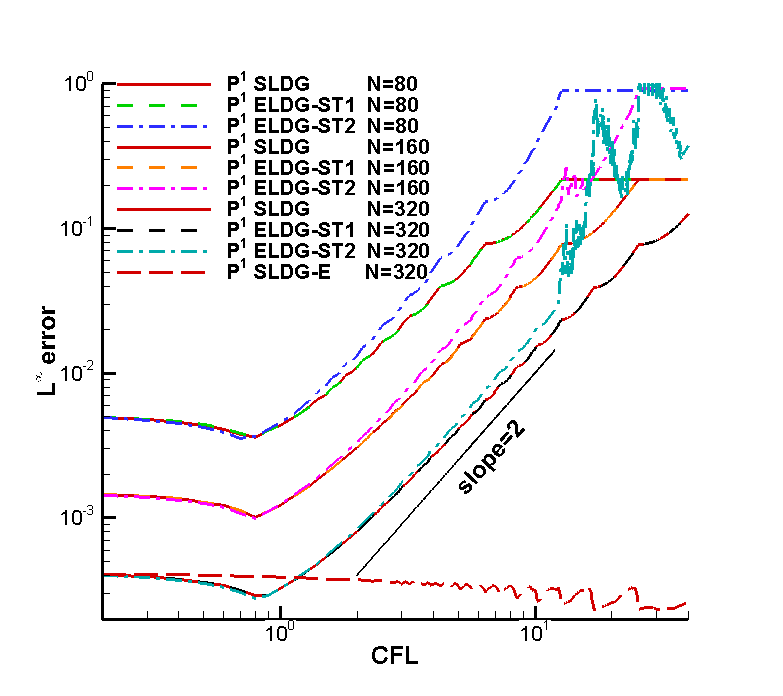}
\includegraphics[height=50mm]{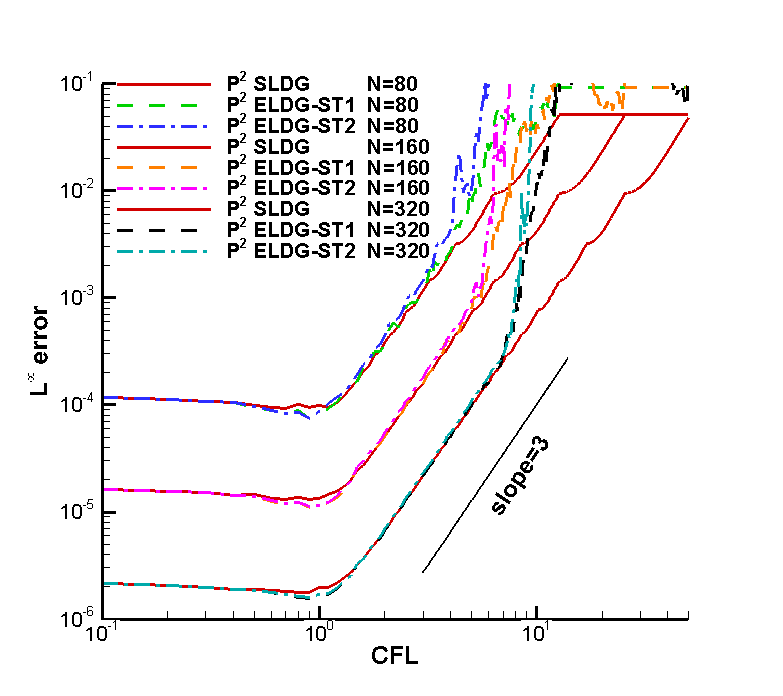}
\caption{The $L^\infty$ error versus $CFL$ of SLDG methods and ELDG methods for 1D transport equation with variable coefficients. $u_t + (\sin(x)u)_x =0$ with the initial condition $u(x,0) = 1$. $T=1$.    $\Delta t=CFL\Delta x$.
$P^1$ SLDG-E means $P^1$ SLDG scheme tracking characteristic lines exactly.
}
\label{sin1d_time}
\end{figure}



\end{example}

\subsection{2D linear transport problems}

\begin{example} \label{example:rotation}
(Rigid body rotation.) Consider
\begin{equation}
u_t -(yu)_x +(xu)_y=0,\ (x,y)\in[-\pi,\pi]^2 .
\end{equation}
The initial condition is set to be the following smooth cosine bell (with $C^5$ smoothness),
\begin{equation}
u(x,y,0) =
\begin{cases}
r_0^b \cos^6 \left(  \frac{r^b}{2r_0^b} \pi \right), & \text{if}\  r^b <r_0^b,\\
0,  & \text{ otherwise},
\end{cases}
\label{cosine_bell}
\end{equation}
 where $r_0^b = 0.3\pi$, and $r^b=\sqrt{  (x-x_0^b)^2 + (y-y_0^b)^2 }$ denotes the distance between $(x,y)$ and the center of the cosine bell $(x_0^b, y_0^b) = ( 0.3\pi ,0 )$.
 First of all, we present the spatial accuracies of ELDG, SLDG and RKDG for solving this problem up to $T=2\pi$ in Table \ref{table:rotation}; the expected $k+1$th order  of convergence   is observed for these schemes with $P^k$ polynomial space. Then, we study numerical stabilities of ELDG and SLDG methods.
In Figure \ref{rotation_cfl}, we present the plots of $L^\infty$ error versus $CFL$ of ELDG and SLDG schemes with different meshes. A few observations can be made:  (1) When $CFL$ is around and below order $1$, both schemes have similar performance in error magnitude and order of convergence. Notice that this time stepping size is well above the stability constraint of $1/(2k+1)$ for RKDG. (2) When $CFL$ is relatively large but smaller than the stability constraint of ELDG, the temporal error starts to kick in 2nd and 3rd order temporal convergence order is shown. (3) Maximum $CFL$s of $P^2$ ELDG-ST1 using $N=40,80,160$ are  around $9$, $13$, $18$. The increasing rate is around $1.4$. Maximum CFLs of $P^2$ ELDG-ST2 using $N=40,80,160$ are  around $8$, $11.5$, $16.5$. The increasing rate is around $1.4$. The increasing ratio of upper bounds of $CFL$ is around $\sqrt{2}$, which coincides with $\Delta t\sim \sqrt{\Delta x}$ as in Remark ~\ref{rem:2.5}.
Similar observations can be made for the $P^1$ case.

\begin{table}[!ht]
\footnotesize
\caption{Rigid body rotation. $u_t -(yu)_x +(xu)_y=0$ with the smooth cosine bell. $T=2\pi.$
We use $CFL=0.3$ and $CFL=0.18$ for all $P^1$ and $P^2$ schemes, respectively.
  }
\centering
\begin{tabular}{| c | cc  | cc| cc| cc|  }
\hline
 Mesh  &{$L^\infty$ error} & Order    &  {$L^\infty$ error} & Order &  {$L^\infty$ error} & Order &  {$L^\infty$ error} & Order \\
 \hline
  & \multicolumn{2}{c|}{$P^1$ RKDG} & \multicolumn{2}{c|}{$P^1$ SLDG}& \multicolumn{2}{c|}{$P^1$ ELDG-ST1} & \multicolumn{2}{c|}{$P^1$ ELDG-ST2}
  \\ \hline
    $20^2$ &     5.40E-01 &   -- &    5.53E-01 &    --       &     5.41E-01 &    -- &     5.41E-01 &   --\\
    $40^2$ &     2.47E-01 &     1.13  &   2.59E-01 &     1.09 &     2.47E-01 &     1.13 & 2.47E-01 &     1.13\\
   $80^2$ &     6.17E-02 &     2.00 &    6.64E-02 &     1.96  &    6.17E-02 &     2.00 & 6.17E-02 &     2.00\\
   $160^2$ &     1.03E-02 &     2.58 &    1.11E-02 &     2.58 &     1.03E-02 &     2.58 &      1.03E-02 &     2.58\\
\hline
  & \multicolumn{2}{c|}{$P^2$ RKDG} & \multicolumn{2}{c|}{$P^2$ SLDG-QC}& \multicolumn{2}{c|}{$P^2$ ELDG-ST1}& \multicolumn{2}{c|}{$P^2$ ELDG-ST2}
  \\ \hline
    $20^2$ &    1.49E-01 &  --       &   1.54E-01 &     --        &   1.49E-01 & -- &     1.49E-01 & --\\
    $40^2$ &    1.39E-02 &     3.42  &    1.48E-02 &     3.39  &    1.39E-02 &     3.42 &     1.39E-02 &     3.42\\
   $80^2$ &    1.61E-03 &     3.11   &   1.65E-03 &     3.16   &   1.61E-03 &     3.11 &     1.61E-03 &     3.11 \\
   $160^2$ &     2.18E-04 &     2.89  &   2.23E-04 &     2.89  &    2.18E-04 &     2.89 &     2.18E-04 &     2.89 \\
\hline
\end{tabular}
\label{table:rotation}
\end{table}

\begin{figure}[!ht]
\centering
\includegraphics[height=50mm]{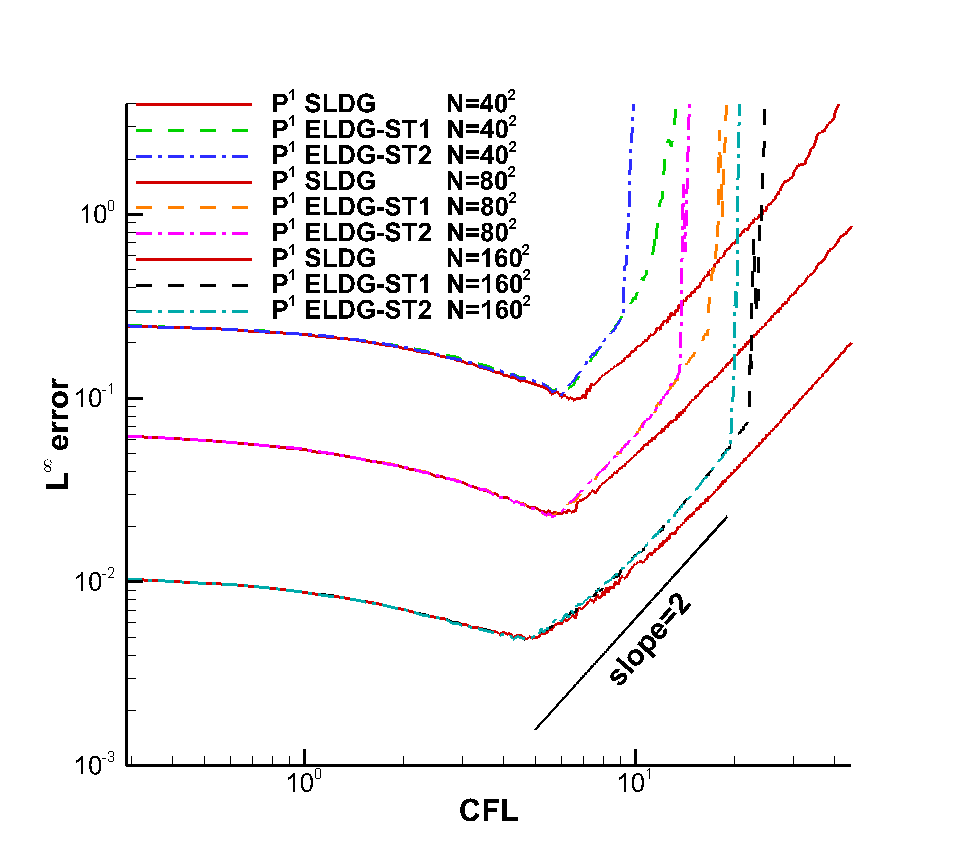}
\includegraphics[height=50mm]{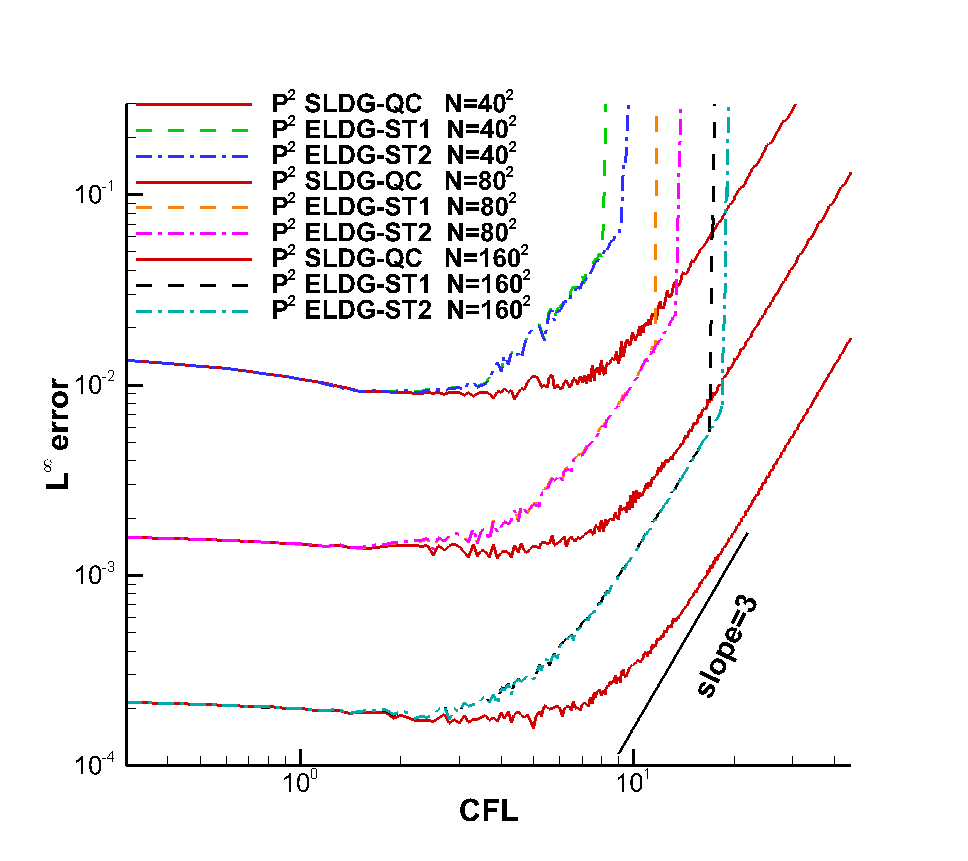}
\caption{
The $L^\infty$ error versus $CFL$ of SLDG schemes and ELDG schemes for the rigid body rotation  with  the smooth cosine bells \eqref{cosine_bell}. 
$T=2\pi$.
 }
\label{rotation_cfl}
\end{figure}

\end{example}

\begin{example}
(Swirling deformation flow.) We consider solving
\begin{equation}
u_t - \left( \cos^2\left(\frac{x}{2} \right) \sin(y) g(t) u \right)_x + \left( \sin(x) \cos^2\left(\frac{y}{2} \right) g(t) u \right)_y =0,\
(x,y)\in[-\pi,\pi]^2,
\end{equation}
with the same initial condition \eqref{cosine_bell},
where $g(t)=\cos\left( \frac{\pi t}{T} \right)\pi$ and $T=1.5$.
As Example \ref{example:rotation}, we also study the spatial error and the numerical stability of the proposed ELDG schemes  in Table \ref{swirling:spatial} and Figure \ref{swirling:temporal}, respectively.
The similar observations as Example \ref{example:rotation} can be made.

\begin{table}[!ht]
\footnotesize
\caption{Swirling deformation flow with  the smooth cosine bells \eqref{cosine_bell}.  $T=1.5$.
We use $CFL=0.3$ and $CFL=0.18$ for all $P^1$ and $P^2$ schemes, respectively.
  }
\centering
\begin{tabular}{| c | cc  | cc| cc| cc|  }
\hline
 Mesh  &{$L^\infty$ error} & Order    &  {$L^\infty$ error} & Order &  {$L^\infty$ error} & Order &  {$L^\infty$ error} & Order \\
 \hline
  & \multicolumn{2}{c|}{$P^1$ RKDG} & \multicolumn{2}{c|}{$P^1$ SLDG}& \multicolumn{2}{c|}{$P^1$ ELDG-ST1} & \multicolumn{2}{c|}{$P^1$ ELDG-ST2}
  \\ \hline
    $20^2$ &     4.00E-01 &   -- &    3.76E-01 &    --       &     3.77E-01 &    -- &    3.76E-01 &   --\\
    $40^2$ &     1.55E-01 &     1.37  &   1.39E-01 &     1.43 &     1.39E-01 &     1.44 & 1.39E-01 &     1.44\\
   $80^2$ &     3.54E-02 &     2.13  &    3.15E-02 &     2.15  &    3.14E-02 &     2.15 & 3.13E-02 &     2.15\\
   $160^2$ &     6.29E-03 &     2.49 &    5.62E-03 &     2.49 &     5.58E-03 &     2.49 &    5.57E-03 &     2.49\\
\hline
  & \multicolumn{2}{c|}{$P^2$ RKDG} & \multicolumn{2}{c|}{$P^2$ SLDG-QC}& \multicolumn{2}{c|}{$P^2$ ELDG-ST1}& \multicolumn{2}{c|}{$P^2$ ELDG-ST2}
  \\ \hline
    $20^2$ &   9.80E-02 &  --       &   9.12E-02 &    --        &   8.97E-02 & -- &      8.92E-02 & --\\
    $40^2$ &    1.33E-02 &     2.88  &   1.13E-02 &     3.02  &    1.04E-02 &     3.11 &     1.04E-02 &     3.10\\
   $80^2$ &    1.79E-03 &     2.89   &   1.58E-03 &     2.84   &   1.47E-03 &     2.82 &     1.47E-03 &     2.82 \\
   $160^2$ &    2.28E-04 &     2.97  &   2.08E-04 &     2.93  &    1.98E-04 &     2.90 &     1.98E-04 &     2.89  \\
\hline
\end{tabular}
\label{swirling:spatial}
\end{table}

\begin{figure}[!ht]
\centering
\includegraphics[height=50mm]{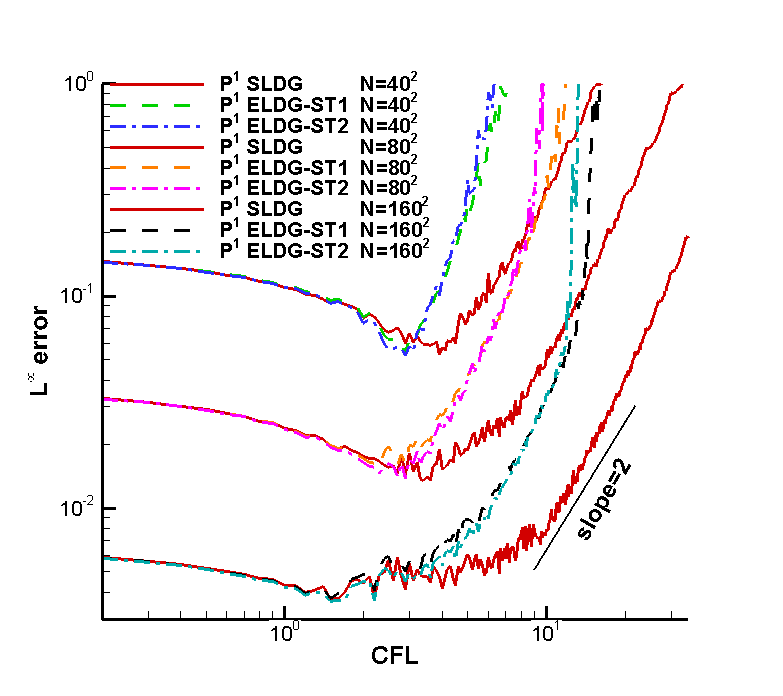}
\includegraphics[height=50mm]{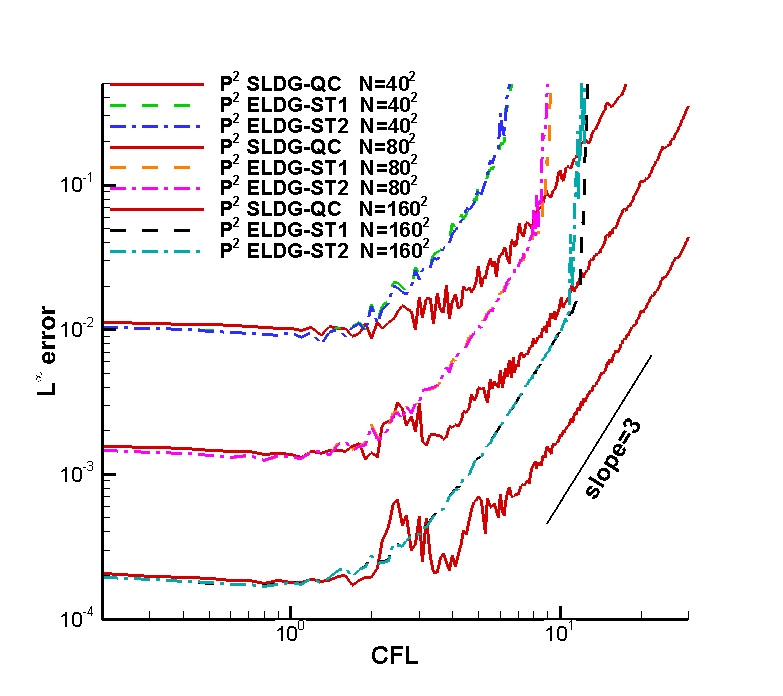}
\caption{
The $L^\infty$ error versus $CFL$ of SLDG methods and ELDG methods for the swirling deformation flow with  the smooth cosine bells \eqref{cosine_bell} with $T=1.5$.
 }
\label{swirling:temporal}
\end{figure}

\end{example}

\subsection{Vlasov-Poisson system}

\begin{example}
(Vlasov-Poisson system: strong Landau damping.)
Consider the strong Landau damping for the
Vlasov-Poisson system \eqref{vlasov} with the initial condition being a perturbed equilibrium
\begin{equation}
f(x,v,t=0) = \frac{ 1 }{\sqrt{2\pi} } ( 1+ \alpha \cos(kx)  ) \exp \left( -\frac{v^2}{2} \right) ,
\label{init}
\end{equation}
with $\alpha = 0.5$ and $k=0.5$ on a computational domain, $[0,4\pi]\times[-2\pi,2\pi]$.
There are several invariants of this problem which should remain constant in time.
These include $L^p$ norms, kinetic energy and entropy:
\begin{itemize}
  \item $L^p$ norm, $1\leq p \leq \infty$:
      \begin{equation}
        \| f \|_p = \left(  \int_v \int_x |f(x,v,t)|^p\,dxdv \right)^{\frac{1}{p}},
      \end{equation}
  \item Energy:
      \begin{equation}
        \text{Energy} = \int_v \int_x f(x,v,t) v^2 dx dv + \int_x E^2 (x,t)dx,
      \end{equation}
  \item Entropy:
     \begin{equation}
      \text{Entropy} = \int_v \int_x f(x,v,t) \log (f(x,v,t) ) dxdv.
     \end{equation}
\end{itemize}
This is a classical  problem that has been numerically investigated by several authors, e.g. see \cite{xiong2014high,zhu2016h,huang2016semi,cai2018high}.

We first test the spatial accuracy of ELDG with the third order temporal scheme for this problem and report the results in Table \ref{VP_spatial}.
The  time reversibility of the Vlasov-Poisson system \cite{degond2004modeling} is used to test the order of convergence.
In Table \ref{VP_spatial}, we show the $L^1$errors and the corresponding orders  of convergence for $P^k$ ELDG and SLDG, $k=1,2$ with $CFL=0.1$.
We observe the expected orders of convergence of ELDG and SLDG.

We then test the numerical stability of ELDG schemes with different meshes for this problem integrated to $T=5$.
Figure \ref{VP_time} reports $L^\infty$ errors versus $CFL$ of solutions of ELDG schemes as well as the SLDG scheme.
From this Figure, we find the expected orders of convergence of the temporal schemes; we also find that the scheme can allow for as large as $CFL=50$;
we observe that the results of ELDG are very close to those of SLDG.


We next study the performances of ELDG for conserving invariants of this problem.
The parameters of the tests are set as follows: we use a mesh of $160\times160$ cells and $CFL=10$. 
For mass conservation, we observed that the mass deviation of ELDG schemes is   around $-4\times 10^{-9}$  due to the domain cut-off in the velocity space; we omit this result.
Figure \ref{VP_norms} shows time evolutions of the relative deviation of  $L^2$ norms of the solution as well as the discrete kinetic energy and entropy.
We make the observations for this Figure:
$P^2$ ELDG performs better than $P^1$ ELDG for conserving $L^2$ norm, as SLDG schemes;
for conserving energy,  ELDG is worse than SLDG;
for conserving entropy, ELDG does a better job than SLDG.

Finally, we study  ELDG schemes for this problem for a long-time simulation.
We present the plots of solutions of ELDG schemes at $T=40$ in the middle and right panels of Figure \ref{VP_time}.
We observe that $P^2$ ELDG performs much better than $P^1$ ELDG for capturing the filamentation structures.
We find that the solutions of both $P^1$ and $P^2$ ELDG are negative around the places where the density is close to vacuum.
Therefore, the positivity-preserving limiter should be added to the current scheme, for which we plan to explore in the future.

\begin{table}[!ht]
\footnotesize
\caption{Strong Landau damping.  $T=0.5$. Use the time reversibility of the VP system. Order of accuracy in space for the SLDG method and the ELDG method. The third order temporal scheme CF3C03 is used for all schemes. We set $CFL=0.1$ so that the spatial error is the dominant error.
  }
\centering
\begin{tabular}{| c | cc  | cc| cc| cc|  }
\hline
 Mesh  &{$L^1$ error} & Order    &  {$L^1$ error} & Order &  {$L^1$ error} & Order &  {$L^1$ error} & Order \\
 \hline
  & \multicolumn{2}{c|}{$P^1$ SLDG} & \multicolumn{2}{c|}{$P^1$ ELDG}& \multicolumn{2}{c|}{$P^2$ SLDG-QC} & \multicolumn{2}{c|}{$P^2$ ELDG}
  \\ \hline
    $32^2$ &    5.88E-04 &   -- &    5.90E-04 &    --       &     3.69E-05 &    -- &    3.25E-05 &  --\\
    $64^2$ &    1.50E-04 &    1.97  &   1.51E-04 &     1.97 &     4.39E-06 &    3.07 & 3.82E-06 &     3.09 \\
   $96^2$ &     6.67E-05 &    1.9  &    6.71E-05 &     1.99  &    1.28E-06 &    3.04 & 1.11E-06 &     3.04 \\
   $128^2$ &     3.76E-05 &    2.00 &    3.78E-05 &     2.00 &     5.37E-07 &    3.02 &    4.66E-07 &     3.03\\
   $160^2$ &     2.41E-05 &    2.00 &    2.42E-05 &     2.00 &   2.74E-07 &    3.02 &     2.38E-07 &     3.02 \\
\hline
\end{tabular}
\label{VP_spatial}
\end{table}

 \begin{figure}[!ht]
\centering
\includegraphics[height=37mm]{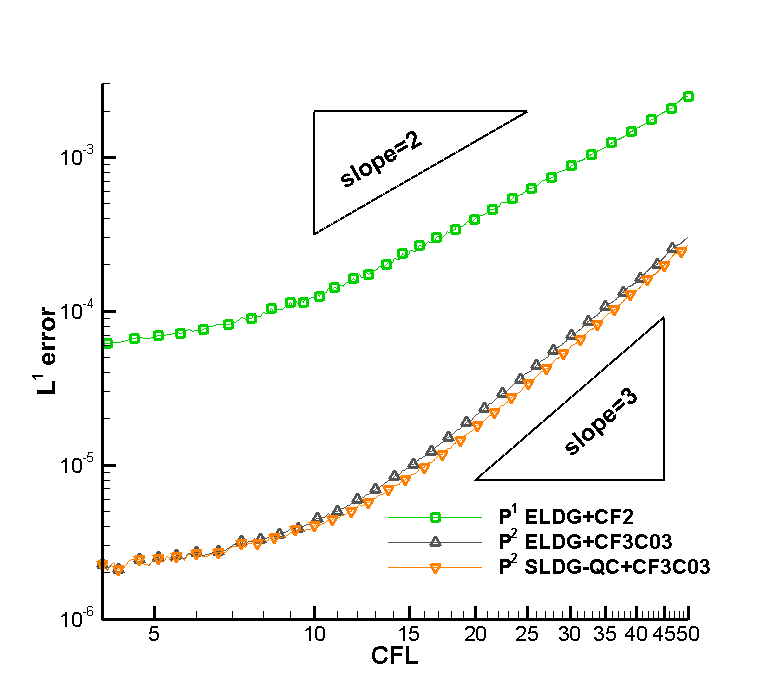}
\includegraphics[height=37mm]{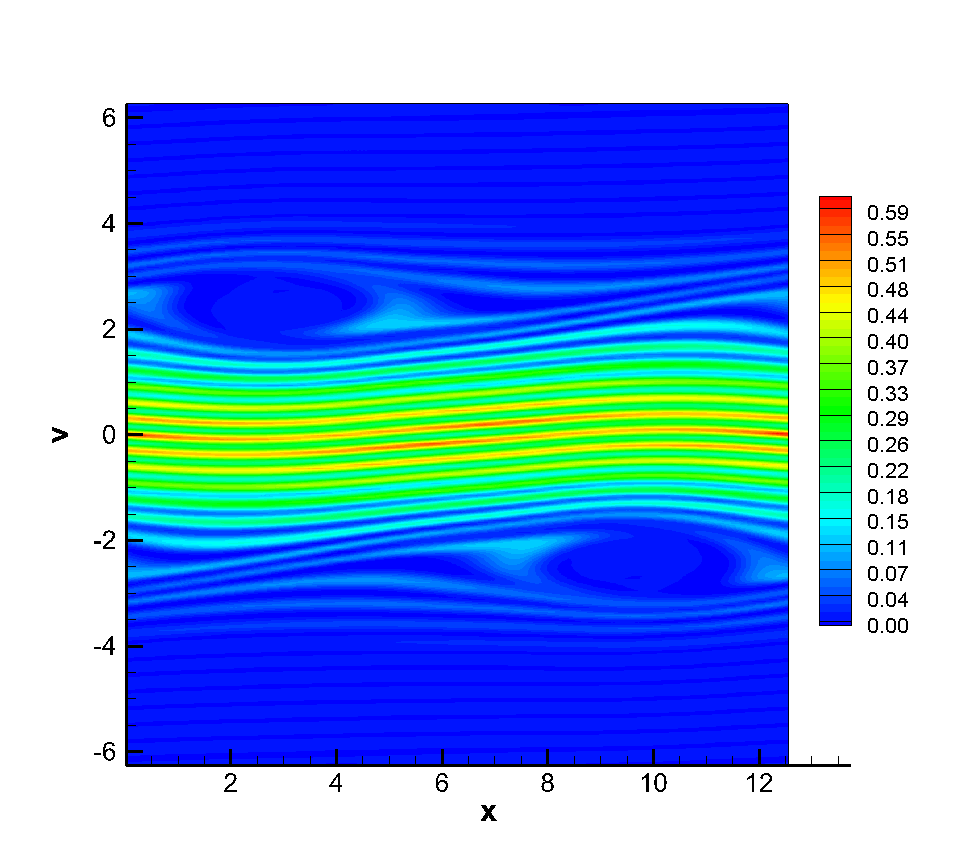}
\includegraphics[height=37mm]{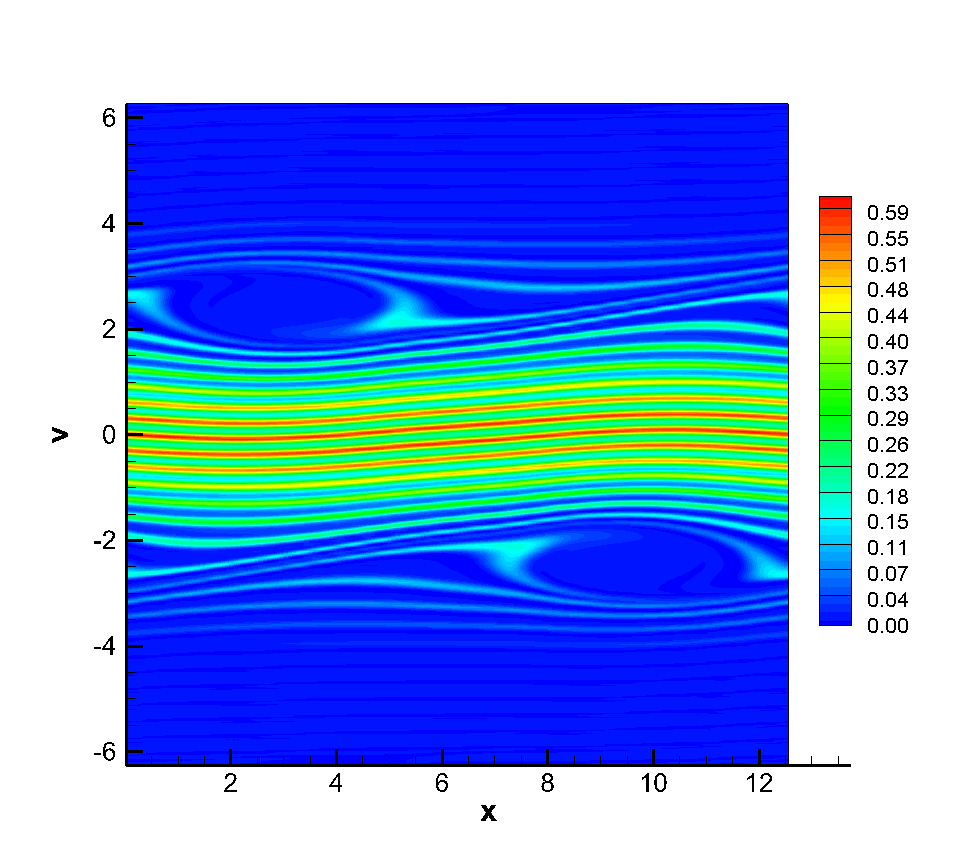}
\caption{Left panel: plots of $L^\infty$ errors versus the $CFL$ number for solving Strong Landau damping at $T=5$.
Temporal order of convergence in  $L^\infty$ norm of ELDG schemes as well as the SLDG scheme coupled with exponential integrators  by comparing numerical solutions with
a reference solution from the corresponding scheme with $CFL = 0.1$.  \protect\\
Middle and right panels: surface plots of the numerical solutions for the strong Landau damping at $T=40$.
We use a mesh of $160\times160$ cells and $CFL=10$. Middle: $P^1$ ELDG+CF2. Right: $P^2$ ELDG+CF3C03.
 }
\label{VP_time}
\end{figure}

\begin{figure}[!ht]
\centering
\includegraphics[height=37mm]{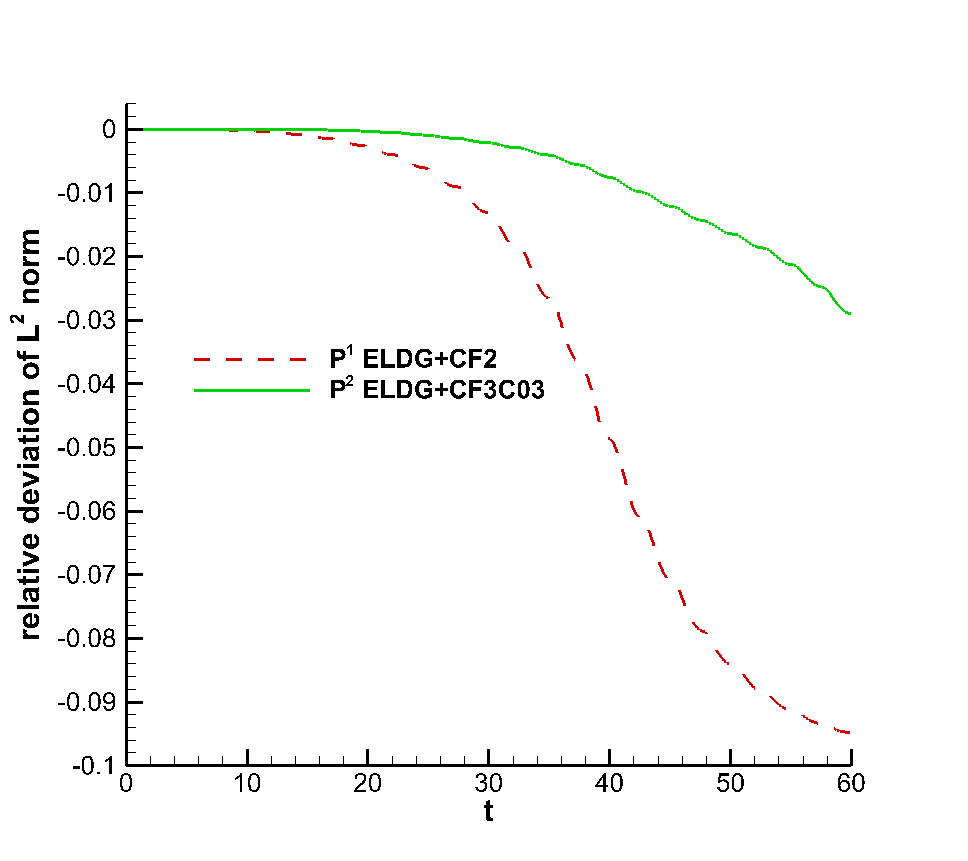}
\includegraphics[height=37mm]{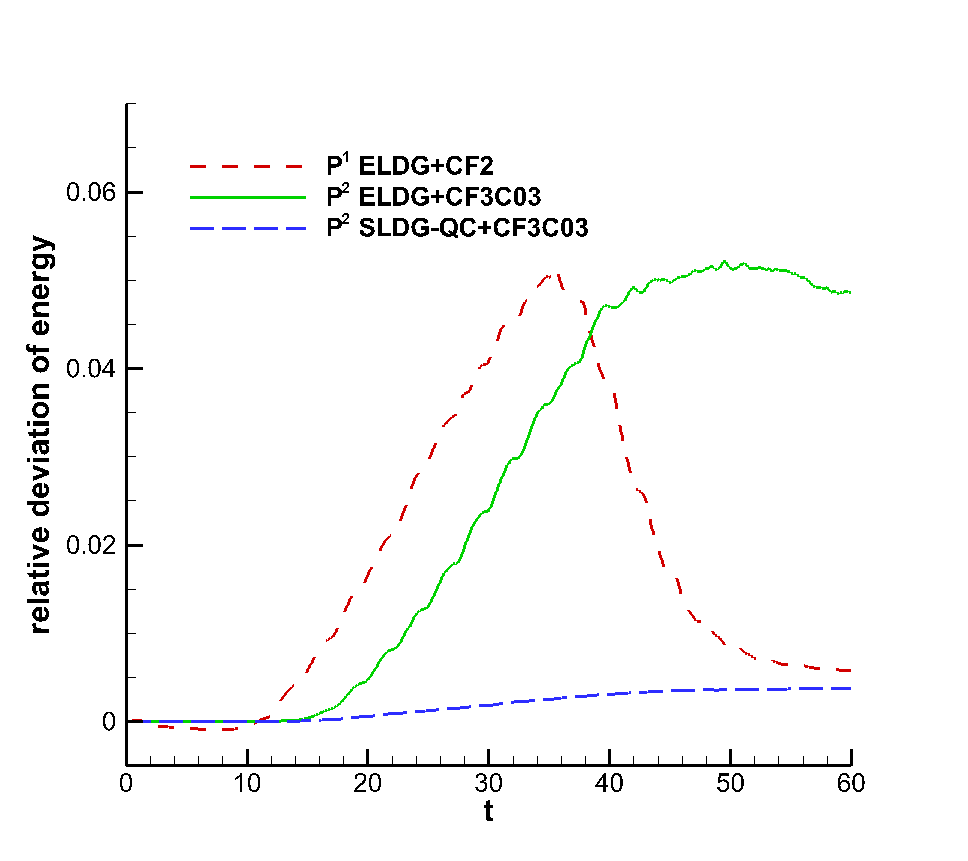}
\includegraphics[height=37mm]{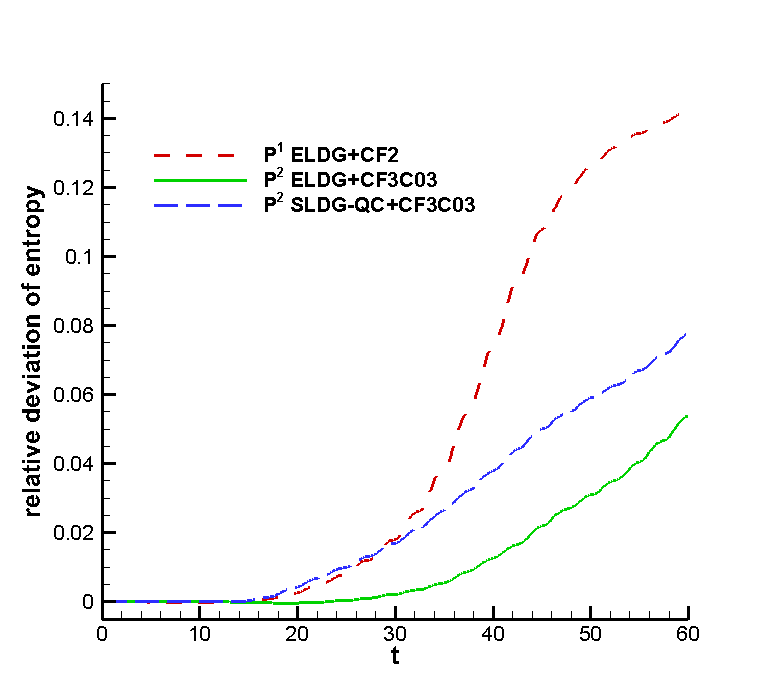}
\caption{Strong Landau damping. Time evolutions of the relative deviation of   $L^2$ (left) norms of the solution as well as the discrete kinetic energy (middle) and entropy (right). We use a mesh of $160\times160$ cells and $CFL=10$ for all simulations. }
\label{VP_norms}
\end{figure}


\end{example}

\subsection{The guiding center Vlasov model}

\begin{example}
(The guiding center Vlasov model: spatial accuracy and convergence test.)
Consider the guiding center Vlasov model on the domain $[0,2\pi]\times[0,2\pi]$ with the initial condition,
$
\rho(x,y,0) = -2 \sin(x) \sin(y)
$
and the periodic boundary condition.
The exact solution stays stationary.
We test the spatial convergence of the proposed ELDG schemes as well as SLDG schemes with the third order temporal scheme, CF3C03,
for solving the guiding center Vlasov model up to time $T=1$ and report the results in Table \ref{guiding_spatial}.
We make the following observations:
(1) we find the expected orders of convergence for $P^k$ ELDG+$P^{k+1}$ LDG, $k=1,2$, in $L^2$ and $L^\infty$ norms;
(2) the results of ELDG schemes are almost the same as those of SLDG schemes.

\begin{table}[!ht]
\footnotesize
\caption{The guiding center Vlasov model on the domain $[0,2\pi]\times[0,2\pi]$ with the initial condition
$
\rho (x,y,0) = -2 \sin(x) \sin(y).
$
Periodic boundary conditions in two directions.
Spatial orders of convergence of $P^{k}$ SLDG(-QC)+$P^{r}$ LDG+CF3C03 and $P^{k}$ ELDG+$P^{r}$ LDG+CF3C03, $k=1,2$, and $r=k+1$.
$T=1$.
$CFL=1$.
  }
\centering
\begin{tabular}{| c | cc    cc| cc  cc|  }
\hline
 Mesh  &{$L^2$ error} & Order    &  {$L^\infty$ error} & Order &  {$L^2$ error} & Order &  {$L^\infty$ error} & Order \\
 \hline
  & \multicolumn{4}{c|}{$P^1$ SLDG} & \multicolumn{4}{c|}{$P^1$ ELDG}
  \\ \hline
    $20^2$ &     1.88E-02 &  -- &     1.06E-01 &    --       &    1.29E-02 & --  &     8.52E-02  &   --\\
    $40^2$ &     4.97E-03 &     1.92 &     3.12E-02 &     1.76 &    3.15E-03 &     2.03 &     2.46E-02 &     1.79  \\
   $60^2$ &     2.24E-03 &     1.97 &     1.44E-02 &     1.90  &    1.36E-03 &     2.07 &     1.14E-02 &     1.90\\
   $80^2$ &     1.27E-03 &     1.95 &     8.27E-03 &     1.93 &     7.71E-04 &     1.98 &     6.52E-03 &     1.93\\
   $100^2$ &    8.17E-04 &     1.99 &     5.34E-03 &     1.96 &     4.94E-04 &     2.00 &     4.22E-03 &     1.95\\
\hline
   & \multicolumn{4}{c|}{$P^2$ SLDG-QC} & \multicolumn{4}{c|}{$P^2$ ELDG}
  \\ \hline
    $20^2$ &     2.77E-03 & -- &     2.06E-02 &    --       &     2.02E-03 & -- &     1.13E-02 &   --\\
    $40^2$ &     3.63E-04 &     2.93 &     4.72E-03 &     2.13 &     2.43E-04 &     3.06 &     2.63E-03 &     2.11\\
   $60^2$ &     1.09E-04 &     2.96 &     2.06E-03 &     2.04  &    7.17E-05 &     3.01 &     1.15E-03 &     2.04\\
   $80^2$ &     4.74E-05 &     2.91 &     1.14E-03 &     2.05 &     2.90E-05 &     3.15 &     6.39E-04 &     2.05\\
   $100^2$ &    2.44E-05 &     2.98 &     7.28E-04 &     2.02 &     1.49E-05 &     2.99 &     4.07E-04 &     2.03\\
\hline
\end{tabular}
\label{guiding_spatial}
\end{table}

\end{example}

\begin{example}
(The guiding center Vlasov model: Kelvin-Helmholtz instability problem.)
We consider the two-dimensional guiding center model problem \eqref{guiding} with the initial condition
\begin{equation}
 \rho_0(x,y) = \sin(y) + 0.015 \cos(kx),
\end{equation}
and periodic boundary condition on the domain $[0,4\pi]\times[0,2\pi]$.
We let $k=0.5$, which will create a Kelvin-Helmholtz instability \cite{shoucri1981two}.

First, we test the temporal convergence of the proposed ELDG schemes with different temporal schemes by computing this problem up to $T=5$.
In particular, we test the proposed second scheme, $P^1$ ELDG+$P^2$ LDG+CF2, and the third order scheme, $P^2$ ELDG+$P^3$ LDG+CF3C03.
In order to minimize the errors for the spatial scheme,  a fixed mesh of $120\times120$ cells is used.
The reference solution is computed by the same scheme with the same mesh but using a small $CFL=0.1$.
We show the plots of $L^1$ errors versus the $CFL$ number of the proposed ELDG schemes for the Kelvin-Helmholtz instability problem at $T=5$ in Figure \ref{KH_temporal}.
We make  a few observations:
(1) we observe expected orders of convergence for all temporal schemes; and $CFL$ of ELDG can be taken to be as large as 50;
(2) by comparing the error magnitude, $P^2$ ELDG+$P^3$ LDG+CF3C03 performs slightly better than $P^2$ SLDG-QC+$P^3$ LDG+CF3C03.
\begin{figure}[!ht]
\centering
\includegraphics[height=45mm]{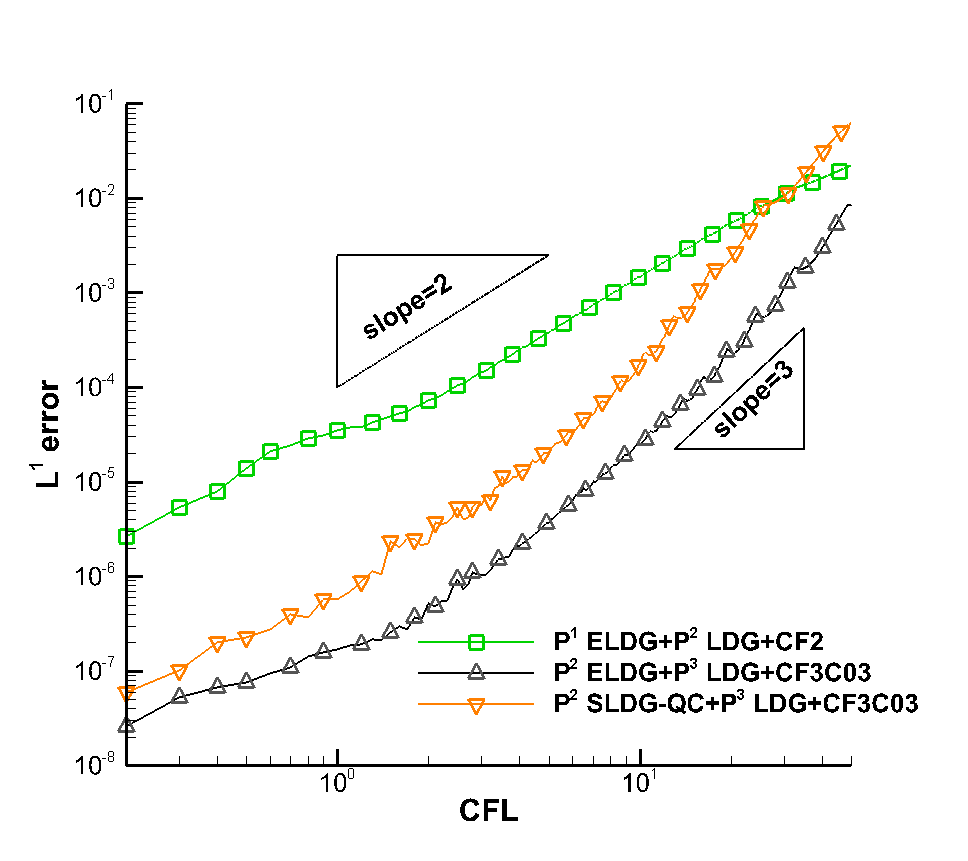}
\caption{ Plots of $L^1$ errors versus the $CFL$ number of the proposed ELDG schemes as well as the SLDG scheme for the Kelvin-Helmholtz instability problem at $T=5$.
Temporal order of convergence of presented schemes by comparing numerical solutions with a reference solution from the corresponding scheme with $CFL=0.1$.
The mesh of $120\times120$ cells is used.}
\label{KH_temporal}
\end{figure}

We then study the quality of the proposed ELDG schemes by tracking relative deviations of some invariants of this problem,
the energy $\| \mathbf{E} \|^2_{L^2} = \int_{\Omega} \mathbf{E} \cdot \mathbf{E} dxdy$ and the enstrophy $\| \rho \|^2_{L^2} = \int_{\Omega} \rho^2 dxdy$.
We study ELDG schemes using a mesh of $100\times100$ cells with $CFL=5$ for solving this problem for a long-time simulation and report the results in Figure \ref{KH_norms}.
We find that $P^2$ ELDG can perform much better than $P^1$ ELDG for conserving both energy and enstrophy.
We find that by comparing SLDG and ELDG with the same polynomial space for conserving both energy and enstrophy, the comparable results can  be observed.
Finally, we  show
surface plots of the numerical solutions for the Kelvin-Helmholtz instability at $T = 40$ in Figure \ref{KH_surface}.
We still observe that the resolution of solutions of ELDG is comparable to that of SLDG.

\begin{figure}[!ht]
\centering
\includegraphics[height=45mm]{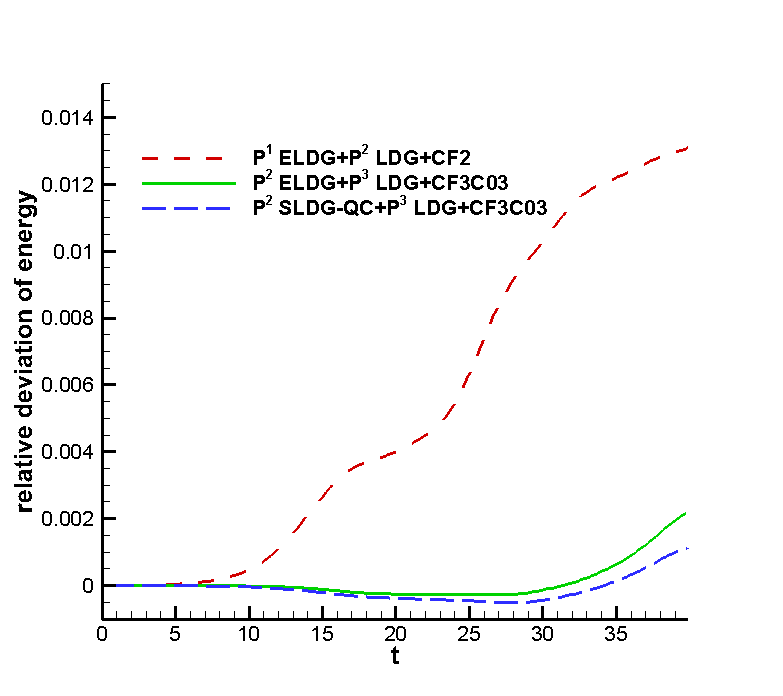}
\includegraphics[height=45mm]{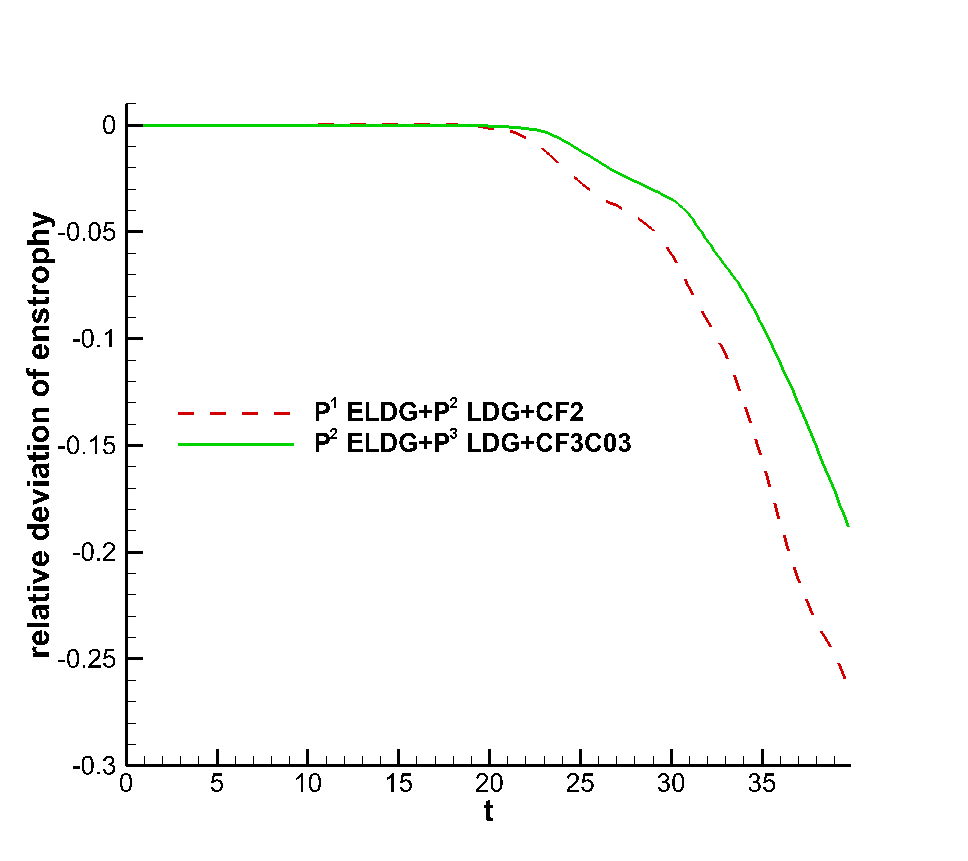}
\caption{Time evolutions of the relative deviation of energy (left) and enstrophy (right) for the proposed ELDG schemes for the Kelvin-Helmholtz instability problem. The mesh of $100\times100$ cells and $CFL=5$ are used. }
\label{KH_norms}
\end{figure}

\begin{figure}[!ht]
\centering
\includegraphics[height=35mm]{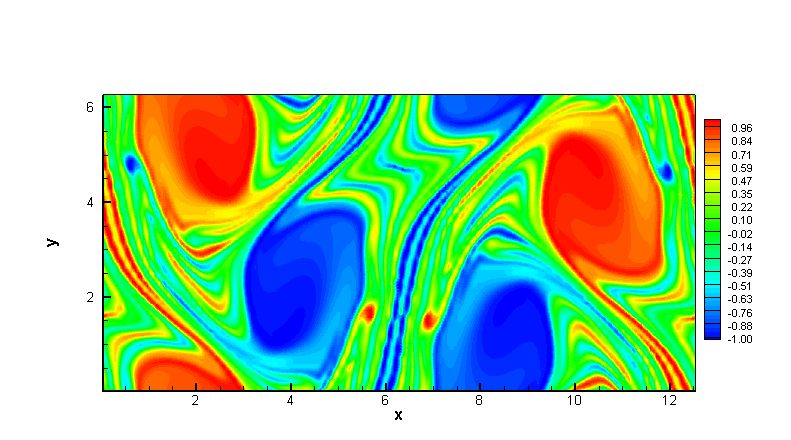}
\includegraphics[height=35mm]{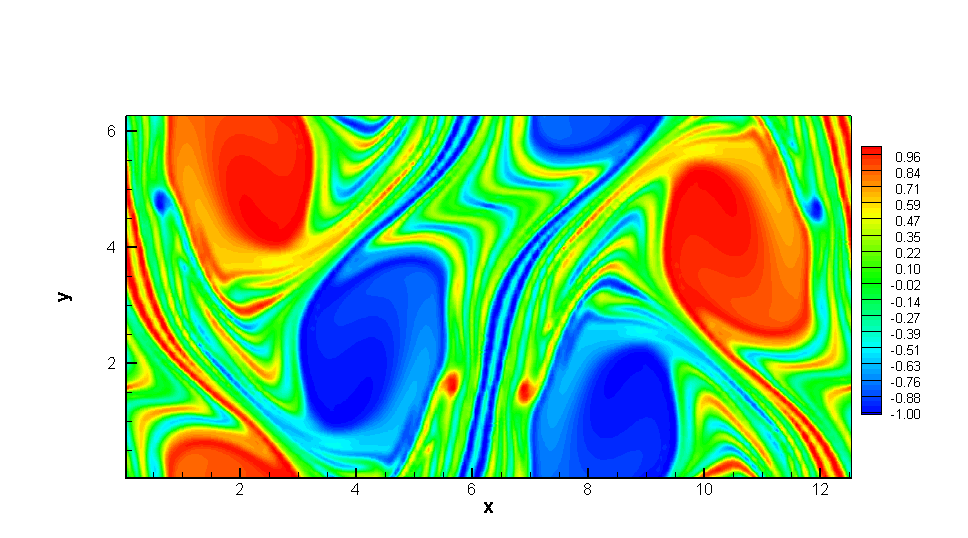}
\caption{Surface plots of the numerical solutions for the Kelvin-Helmholtz instability at $T=40$.
We use a mesh of $100\times100$ cells and $CFL=5$.
Left: $P^2$ SLDG-QC+$P^3$ LDG+CF3C03.
Right: $P^2$ ELDG+$P^3$ LDG+CF3C03.
 }
\label{KH_surface}
\end{figure}

\end{example}

\subsection{The two-dimensional incompressible Euler equations}

\begin{example}
(The incompressible Euler equations: the shear flow problem)
For the double shear layer problem \cite{bell1989second,zhangshu2010}, we solve the 2D incompressible Euler equations \eqref{Euler} in the domain $[0,2\pi]\times[0,2\pi]$, with the periodic boundary conditions and the initial condition given by
\begin{equation}
\omega(x,y,0)
=
\begin{cases}
\delta \cos(x) -\frac{1}{\rho} sech^2\left( \frac{y-\pi/2}{\rho} \right), & \text{if} \ y\leq \pi,\\
\delta \cos(x) + \frac{1}{\rho} sech^2\left( \frac{3\pi/2-y}{\rho} \right), & \text{if}\ y>\pi,
\end{cases}
\end{equation}
where $\delta =0.05$ and $\rho=\pi/15$.

As  time  evolves,  the  solution  quickly  rolls up  with  smaller  and  smaller  spatial scales so on any fixed grid, the full resolution will be lost eventually.
This problem is a classic benchmark for demonstrating the effectiveness of a new scheme so it has been tested for many schemes such as the high order nonsplitting SL WENO scheme \cite{xiong2018high}, the DG method in \cite{liu2000high,zhangshu2010,zhu2017h} and the spectral element method in \cite{fischer2001filter,xu2006stabilization}.
We first show surface plots of numerical solutions for this problem at $T=8$ in Figure \ref{figure:shear1}, where the solution is rolled up in a very small scale.
We find that ELDG schemes could allow for $CFL=5$ for these simulations and the solutions  with larger $CFL=5$ seem to be less dissipative than those with $CFL=1$.
We then study the quality of the ELDG schemes by tracking relative deviations of the energy $\| \mathbf{u} \|^2_{L^2} = \int_{\Omega}\mathbf{u}\cdot \mathbf{u} dxdy$
and the enstrophy $\| \omega \|^2_{L^2} = \int_{\Omega} \omega^2 dxdy$ of this problem and report the results in Figure \ref{figure:shear_norm}.
We observed that higher order $P^2$ ELDG performs much better than the lower order $P^1$ ELDG for conserving energy and enstrophy.

\begin{figure}[!ht]
\centering
\includegraphics[width=50mm]{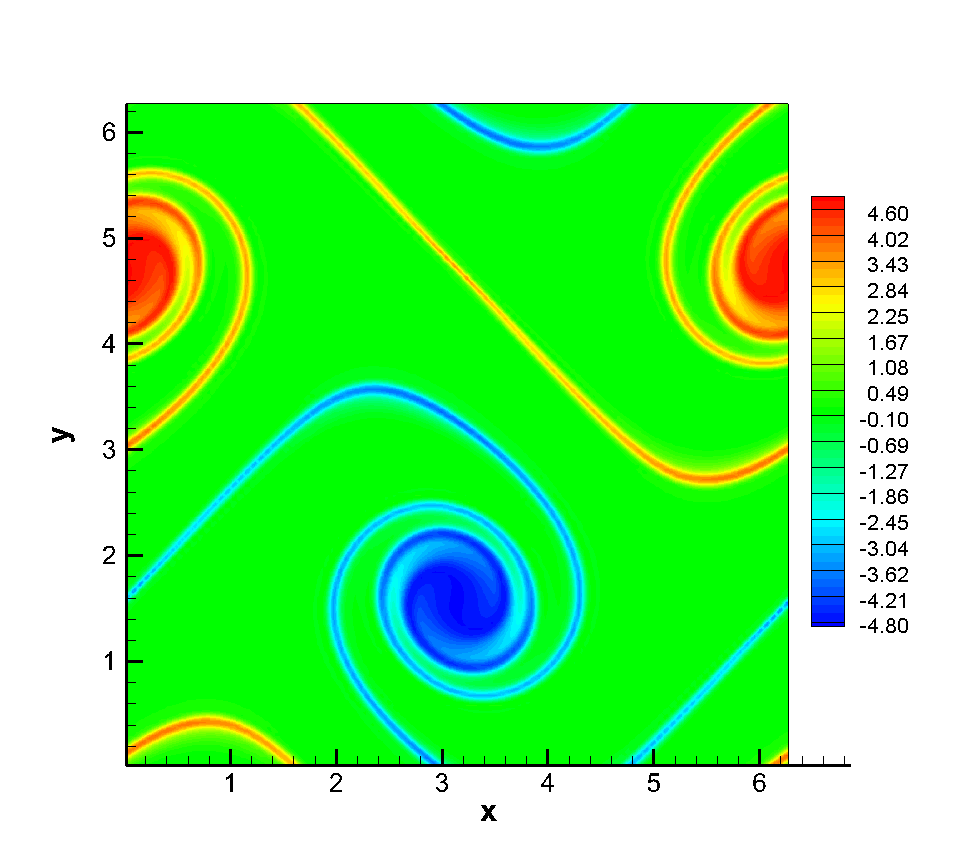}
\includegraphics[width=50mm]{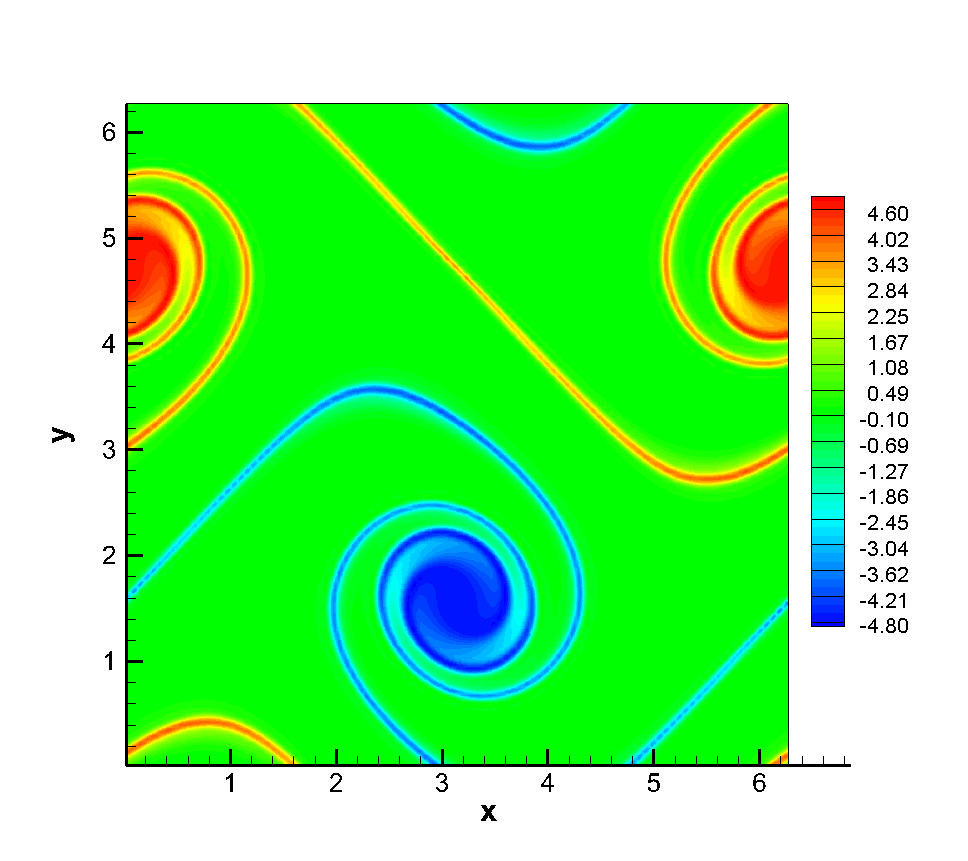}
\caption{Contour plots of the numerical solutions for the shear flow test at $T=8$.
  $P^2$ ELDG +$P^3$ LDG+CF3C03 using $CFL=1$ (left), $CFL=5$ (right).
The mesh of $100\times100$.
}
\label{figure:shear1}
\end{figure}

\begin{figure}[!ht]
\centering
\includegraphics[width=45mm]{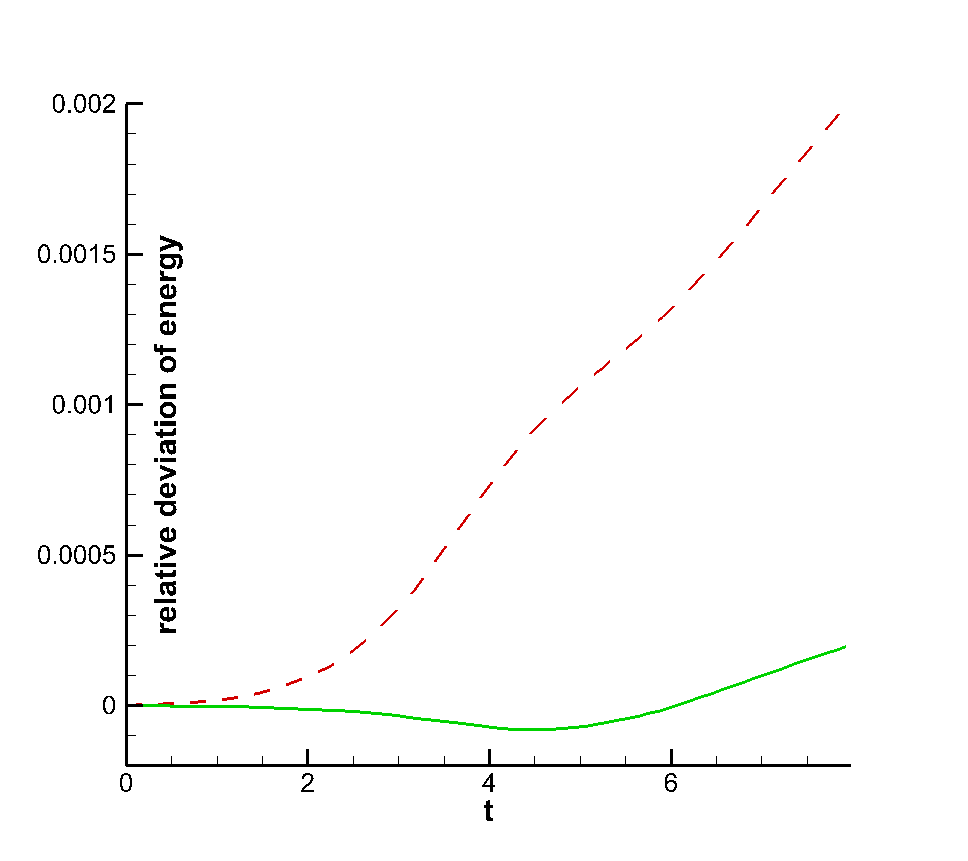}
\includegraphics[width=45mm]{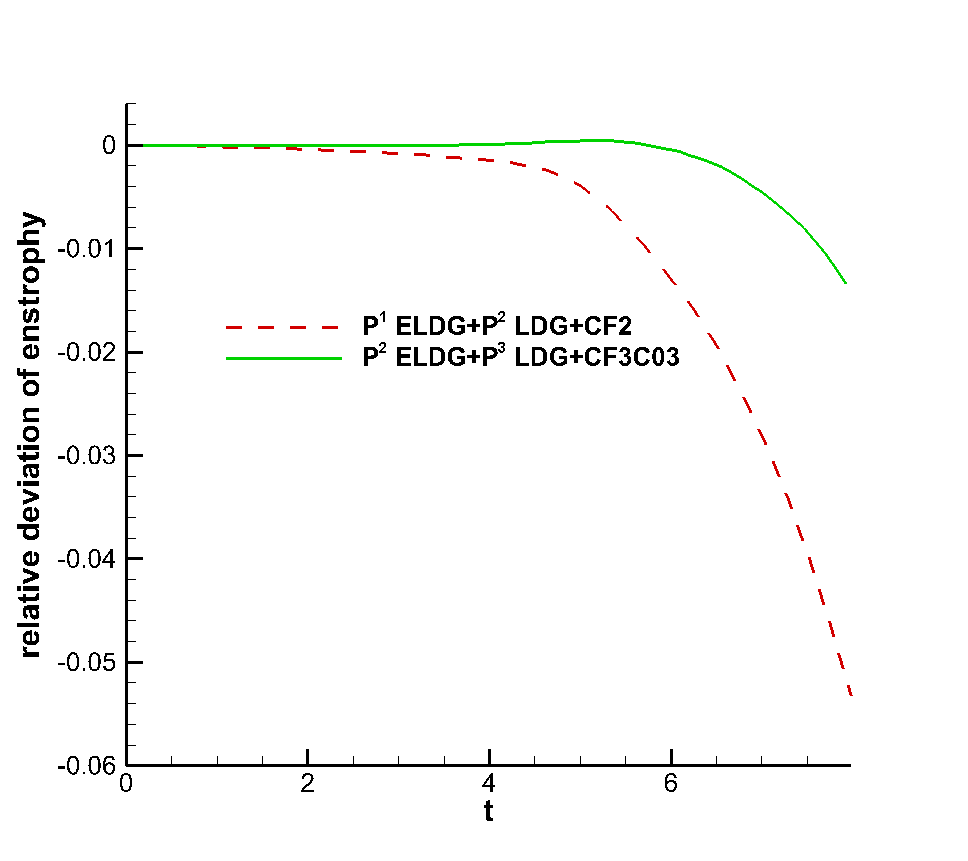}
\caption{Time evolution of the relative deviation of energy and enstrophy for the proposed ELDG schemes for the shear flow test.
Left: ELDG+$P^2$ LDG+CF2.
Right: $P^2$ ELDG +$P^3$ LDG+CF3C03.
We use a mesh of $100\times100$ and $CFL=5$. }
\label{figure:shear_norm}
\end{figure}

\end{example}

\section{Conclusion} \label{section:conclusion}

In this paper, we have developed a new Eulerian-Lagrangian discontinuous Galerkin (DG) method for transport problems. 
The new framework encompasses the semi-Lagrangian DG and Eulerian Runge-Kutta DG in special cases; thus inherits advantages from both approaches in stability under large time stepping sizes, and in mass conservation, compactness and high order accuracy.  
These advantages are numerically verified by extensive numerical tests for linear transport equation and nonlinear dynamics.
 Future works include further theoretic development and application of limiters, developing schemes for nonlinear hyperbolic problems and to unstructured meshes.

\bibliographystyle{siamplain}
\bibliography{references}

\begin{thebibliography}{10}

\bibitem{arnold2002unified}
{\sc D.~Arnold, F.~Brezzi, B.~Cockburn, and L.~Marini}, {\em {Unified analysis
  of discontinuous Galerkin methods for elliptic problems}}, SIAM Journal on
  Numerical Analysis, 39 (2002), pp.~1749--1779.

\bibitem{bell1989second}
{\sc J.~Bell, P.~Colella, and H.~Glaz}, {\em {A second-order projection method
  for the incompressible Navier-Stokes equations}}, Journal of Computational
  Physics, 85 (1989), pp.~257--283.

\bibitem{bosler2019conservative}
{\sc P.~A. Bosler, A.~M. Bradley, and M.~A. Taylor}, {\em {Conservative
  Multimoment Transport along Characteristics for Discontinuous Galerkin
  Methods}}, SIAM Journal on Scientific Computing, 41 (2019), pp.~B870--B902.

\bibitem{cai2019exp}
{\sc X.~Cai, S.~Boscarino, and J.-M. Qiu}, {\em {High Order Semi-Lagrangian
  Discontinuous Galerkin Method Coupled with Runge-Kutta Exponential
  Integrators for Nonlinear Vlasov Dynamics}}, arXiv preprint arXiv:1911.12229,
   (2019).

\bibitem{cai2016high}
{\sc X.~Cai, W.~Guo, and J.-M. Qiu}, {\em {A high order conservative
  semi-Lagrangian discontinuous Galerkin method for two-dimensional transport
  simulations}}, Journal of Scientific Computing, 73 (2017), pp.~514--542.

\bibitem{cai2018high}
{\sc X.~Cai, W.~Guo, and J.-M. Qiu}, {\em {A high order semi-Lagrangian
  discontinuous Galerkin method for Vlasov-Poisson simulations without operator
  splitting}}, Journal of Computational Physics, 354 (2018), pp.~529--551.

\bibitem{castillo2000priori}
{\sc P.~Castillo, B.~Cockburn, I.~Perugia, and D.~Sch{\"o}tzau}, {\em {An a
  priori error analysis of the local discontinuous Galerkin method for elliptic
  problems}}, SIAM Journal on Numerical Analysis, 38 (2000), pp.~1676--1706.

\bibitem{celia1990eulerian}
{\sc M.~Celia, T.~Russell, I.~Herrera, and R.~Ewing}, {\em {An
  Eulerian-Lagrangian localized adjoint method for the advection-diffusion
  equation}}, Advances in Water Resources, 13 (1990), pp.~187--206.

\bibitem{celledoni2009semi}
{\sc E.~Celledoni and B.~K. Kometa}, {\em {Semi-Lagrangian Runge-Kutta
  exponential integrators for convection dominated problems}}, Journal of
  Scientific Computing, 41 (2009), pp.~139--164.

\bibitem{celledoni2003commutator}
{\sc E.~Celledoni, A.~Marthinsen, and B.~Owren}, {\em {Commutator-free Lie
  group methods}}, Future Generation Computer Systems, 19 (2003), pp.~341--352.

\bibitem{ciarlet1988mathematical}
{\sc P.~G. Ciarlet}, {\em {Mathematical Elasticity: Volume I: three-dimensional
  elasticity}}, North-Holland, 1988.

\bibitem{cockburn1989tvb}
{\sc B.~Cockburn and C.-W. Shu}, {\em {TVB Runge-Kutta local projection
  discontinuous Galerkin finite element method for conservation laws II:
  general framework}}, Mathematics of Compututation,  (1989), pp.~411--435.

\bibitem{cockburn1991runge}
{\sc B.~Cockburn and C.-W. Shu}, {\em {The Runge-Kutta local
  projection-discontinuous-Galerkin finite element method for scalar
  conservation laws}}, ESAIM: Mathematical Modelling and Numerical Analysis, 25
  (1991), pp.~337--361.

\bibitem{cockburn1998local}
{\sc B.~Cockburn and C.-W. Shu}, {\em {The local discontinuous Galerkin method
  for time-dependent convection-diffusion systems}}, SIAM Journal on Numerical
  Analysis, 35 (1998), pp.~2440--2463.

\bibitem{cockburn2001runge}
{\sc B.~Cockburn and C.-W. Shu}, {\em {Runge--Kutta discontinuous Galerkin
  methods for convection-dominated problems}}, Journal of Scientific Computing,
  16 (2001), pp.~173--261.

\bibitem{crouseilles2009conservative}
{\sc N.~Crouseilles, M.~Mehrenberger, and E.~Sonnendr{\"u}cker}, {\em
  {Conservative semi-Lagrangian schemes for Vlasov equations}}, Journal of
  Computational Physics, 229 (2010), pp.~1927--1953.

\bibitem{degond2004modeling}
{\sc P.~Degond, L.~Pareschi, and G.~Russo}, {\em Modeling and Computational
  Methods for Kinetic Equations}, Springer, 2004.

\bibitem{fischer2001filter}
{\sc P.~Fischer and J.~Mullen}, {\em Filter-based stabilization of spectral
  element methods}, Comptes Rendus de l'Acad{\'e}mie des Sciences-Series
  I-Mathematics, 332 (2001), pp.~265--270.

\bibitem{Guo2013discontinuous}
{\sc W.~Guo, R.~Nair, and J.-M. Qiu}, {\em A conservative semi-{L}agrangian
  discontinuous {G}alerkin scheme on the cubed-sphere}, Monthly Weather Review,
  142 (2013), pp.~457--475.

\bibitem{huang2017eulerian}
{\sc C.-S. Huang and T.~Arbogast}, {\em {An Eulerian--Lagrangian Weighted
  Essentially Nonoscillatory scheme for Nonlinear Conservation Laws}},
  Numerical Methods for Partial Differential Equations, 33 (2017),
  pp.~651--680.

\bibitem{huang2018implicit}
{\sc C.-S. Huang and T.~Arbogast}, {\em {An Implicit Eulerian--Lagrangian WENO3
  Scheme for Nonlinear Conservation Laws}}, Journal of Scientific Computing, 77
  (2018), pp.~1084--1114.

\bibitem{huang2016semi}
{\sc C.-S. Huang, T.~Arbogast, and C.-H. Hung}, {\em {A semi-Lagrangian finite
  difference WENO scheme for scalar nonlinear conservation laws}}, Journal of
  Computational Physics, 322 (2016), pp.~559--585.

\bibitem{klingenberg2017arbitrary}
{\sc C.~Klingenberg, G.~Schn{\"u}cke, and Y.~Xia}, {\em {Arbitrary
  Lagrangian-Eulerian discontinuous Galerkin method for conservation laws:
  analysis and application in one dimension}}, Mathematics of Computation, 86
  (2017), pp.~1203--1232.

\bibitem{liu2000high}
{\sc J.-G. Liu and C.-W. Shu}, {\em {A high-order discontinuous Galerkin method
  for 2D incompressible flows}}, Journal of Computational Physics, 160 (2000),
  pp.~577--596.

\bibitem{luo2019quasi}
{\sc D.~Luo, W.~Huang, and J.~Qiu}, {\em {A quasi-Lagrangian moving mesh
  discontinuous Galerkin method for hyperbolic conservation laws}}, Journal of
  Computational Physics, 396 (2019), pp.~544--578.

\bibitem{persson2009discontinuous}
{\sc P.-O. Persson, J.~Bonet, and J.~Peraire}, {\em {Discontinuous Galerkin
  solution of the Navier--Stokes equations on deformable domains}}, Computer
  Methods in Applied Mechanics and Engineering, 198 (2009), pp.~1585--1595.

\bibitem{qiu2011positivity}
{\sc J.-M. Qiu and C.-W. Shu}, {\em {Positivity preserving semi-Lagrangian
  discontinuous Galerkin formulation: Theoretical analysis and application to
  the Vlasov--Poisson system}}, Journal of Computational Physics, 230 (2011),
  pp.~8386--8409.

\bibitem{James}
{\sc J.~Rossmanith and D.~Seal}, {\em {A positivity-preserving high-order
  semi-Lagrangian discontinuous Galerkin scheme for the Vlasov-Poisson
  equations}}, Journal of Computational Physics, 230 (2011), pp.~6203--6232.

\bibitem{russell2002overview}
{\sc T.~F. Russell and M.~A. Celia}, {\em {An overview of research on
  Eulerian--Lagrangian localized adjoint methods (ELLAM)}}, Advances in Water
  resources, 25 (2002), pp.~1215--1231.

\bibitem{shoucri1981two}
{\sc M.~M. Shoucri}, {\em A two-level implicit scheme for the numerical
  solution of the linearized vorticity equation}, International Journal for
  Numerical Methods in Engineering, 17 (1981), pp.~1525--1538.

\bibitem{shu1988efficient}
{\sc C.-W. Shu and S.~Osher}, {\em {Efficient implementation of essentially
  non-oscillatory shock-capturing schemes}}, Journal of Computational Physics,
  77 (1988), pp.~439--471.

\bibitem{wang1999family}
{\sc H.~Wang, R.~Ewing, G.~Qin, S.~Lyons, M.~Al-Lawatia, and S.~Man}, {\em {A
  family of Eulerian-Lagrangian localized adjoint methods for multi-dimensional
  advection-reaction equations}}, Journal of Computational Physics, 152 (1999),
  pp.~120--163.

\bibitem{xiong2014high}
{\sc T.~Xiong, J.-M. Qiu, Z.~Xu, and A.~Christlieb}, {\em {High order maximum
  principle preserving semi-Lagrangian finite difference WENO schemes for the
  Vlasov equation}}, Journal of Computational Physics, 273 (2014),
  pp.~618--639.

\bibitem{xiong2018high}
{\sc T.~Xiong, G.~Russo, and J.-M. Qiu}, {\em High order multi-dimensional
  characteristics tracing for the incompressible euler equation and the
  guiding-center vlasov equation}, Journal of Scientific Computing, 77 (2018),
  pp.~263--282.

\bibitem{xu2006stabilization}
{\sc C.~Xu}, {\em Stabilization methods for spectral element computations of
  incompressible flows}, Journal of Scientific Computing, 27 (2006),
  pp.~495--505.

\bibitem{zhangshu2010}
{\sc X.~Zhang and C.-W. Shu}, {\em {On maximum-principle-satisfying high order
  schemes for scalar conservation laws}}, Journal of Computational Physics, 229
  (2010), pp.~3091--3120.

\bibitem{zhou2019stability}
{\sc L.~Zhou, Y.~Xia, and C.-W. Shu}, {\em {Stability analysis and error
  estimates of arbitrary Lagrangian-Eulerian discontinuous Galerkin method
  coupled with Runge-Kutta time-marching for linear conservation laws}}, ESAIM:
  Mathematical Modelling and Numerical Analysis, 53 (2019), pp.~105--144.

\bibitem{zhu2016h}
{\sc H.~Zhu, J.~Qiu, and J.-M. Qiu}, {\em {An h-adaptive RKDG method for the
  Vlasov--Poisson system}}, Journal of Scientific Computing, 69 (2016),
  pp.~1346--1365.

\bibitem{zhu2017h}
{\sc H.~Zhu, J.~Qiu, and J.-M. Qiu}, {\em {An h-Adaptive RKDG Method for the
  Two-Dimensional Incompressible {E}uler Equations and the Guiding Center
  {V}lasov Model}}, Journal of Scientific Computing, 73 (2017), pp.~1316--1337.

\end{thebibliography}
\end{document}